\RequirePackage{etoolbox}
\csdef{input@path}{{style/}{graphics/}}
\documentclass[MSNbibl,number,citesort,dvips]{arxbj}
\makeatletter
   \@ifpackageloaded{graphicx}{}{\usepackage{graphicx}}
\makeatother

\aid{0}
\volume{21}
\issue{1}
\pubyear{2015}
\firstpage{574}
\lastpage{603}
\doi{10.3150/13-BEJ580} 

\makeatletter

\newtheorem{Theorem}{Theorem}
\newtheorem{Corollary}{Corollary}
\newtheorem{Lemma}{Lemma}
\newremark{Preparation}{Preparation}

\newcommand{\rright}{\right}
\newcommand{\lleft}{\left}
\def\diag{\operatorname{diag}}

\def\N{\mathrm{N}}

\def\tr{\operatorname{tr}}
\def\T{\mathrm{T}}
\renewcommand{\emptyset}{\varnothing}

\def\sfrac#1#2{#1/#2}
\def\vfrac#1#2{(#1)/#2}

\def\sklfrac#1#2{(#1/#2)}

\makeatother

\begin{document}
\begin{frontmatter}

\title{Improved minimax estimation of a multivariate normal mean under
heteroscedasticity}
\runtitle{Minimax estimation under heteroscedasticity}

\begin{aug}
\author{\inits{Z.}\fnms{Zhiqiang}~\snm{Tan}\corref{}\ead[label=e1]{ztan@stat.rutgers.edu}}
\address{Department of Statistics, Rutgers University, 110
Frelinghuysen Road,
Piscataway, NJ 08854, USA. \printead{e1}}
\end{aug}

\received{\smonth{5} \syear{2013}}
\revised{\smonth{11} \syear{2013}}

%
\begin{abstract}
Consider the problem of estimating a multivariate normal mean with a
known variance matrix, which is not necessarily proportional to the
identity matrix. The coordinates are shrunk directly in proportion to
their variances in Efron and Morris' (\textit{J. Amer. Statist. Assoc.}
\textbf{68} (1973) 117--130) empirical Bayes approach,
whereas inversely in proportion to their variances in Berger's (\textit
{Ann. Statist.}
\textbf{4} (1976) 223--226)
minimax estimators. We propose a new minimax estimator, by
approximately minimizing the Bayes risk with a normal prior among a
class of minimax estimators where the shrinkage direction is open to
specification and the shrinkage magnitude is determined to achieve
minimaxity. The proposed estimator has an interesting simple form such
that one group of coordinates are shrunk in the direction of Berger's
estimator and the remaining coordinates are shrunk in the direction of
the Bayes rule.
Moreover, the proposed estimator is scale adaptive: it can achieve
close to the minimum Bayes risk simultaneously over a scale class of
normal priors (including the specified prior) and achieve close to the
minimax linear risk over a corresponding scale class of hyper-rectangles.
For various scenarios in our numerical study, the proposed estimators
with extreme priors yield more substantial risk reduction than existing
minimax estimators.
\end{abstract}

%

\end{frontmatter}

\section{Introduction}

A fundamental statistical problem is shrinkage estimation of a
multivariate normal mean. See, for example, the February 2012 issue of
\textit{Statistical Science} for a broad range of theory, methods, and
applications.
Let $X=(X_1,\ldots, X_p)^\T$ be multivariate normal with \textit
{unknown} mean vector $\theta=(\theta_1,\ldots,\theta_p)^\T$ and
\textit{known} variance matrix $\Sigma$. Consider the problem of
estimating $\theta$ by an estimator $\delta=\delta(X)$ under the loss
$L(\delta, \theta) = (\delta-\theta)^\T Q (\delta-\theta)$,
where $Q$ is a \textit{known} positive definite, symmetric matrix. The
risk of $\delta$ is $R(\delta,\theta)=E_\theta\{ L(\delta,\theta
)\}$.
The general problem can be transformed into a canonical form such that
$\Sigma$ is diagonal and $Q=I$, the identity matrix (e.g., Lehmann and
Casella \cite{LehCas98}, Problem 5.5.11). For simplicity, assume except in
Section~\ref{sec3.2} that $\Sigma$ is $D=\diag(d_1, \ldots, d_p) $ and
$L(\delta,\theta)=\|\delta- \theta\|^2$, where $\|x \|^2 = x^{\T} x$
for a column vector $x$.
The letter $D$ is substituted for $\Sigma$ to emphasize that it is diagonal.

For this problem, we aim to develop shrinkage estimators that are both
minimax and capable of effective risk reduction over the usual
estimator $\delta_0=X$ even in the heteroscedastic case (i.e.,
$d_1,\ldots,d_p$ are not equal). An estimator of $\theta$ is minimax
if and only if, \textit{regardless of} $\theta\in\mathbb R^p$, its
risk is always no greater than $\sum_{j=1}^p d_j$, the risk of $\delta
_0$. For $p\ge3$, minimax estimators different from and hence
dominating $\delta_0$ are first discovered in the homoscedastic case
where $D=\sigma^2 I$ (i.e., $d_1=\cdots=d_p=\sigma^2$).
James and Stein \cite{JamSte61} showed that $\delta_c^{\mathrm{JS}} = (1-c \sigma^2 /\|X
\|^2 ) X$
is minimax provided $0 \le c \le2(p-2)$.
Stein \cite{Ste62} suggested the positive-part estimator
$\delta_c^{\mathrm{JS}+} =
(1-c \sigma^2/\|X \|^2)_+ X$, which dominates $\delta_c^{\mathrm{JS}}$.
Throughout, $a_+ =\max(0, a)$.
Shrinkage estimation has since been developed into a general
methodology with various approaches, including empirical Bayes (Efron
and Morris \cite{EfrMor73}; Morris \cite{Mor83}) and hierarchical
Bayes (Strawderman \cite{Str71};
Berger and Robert \cite{BerRob90}). While these approaches are prescriptive for
constructing shrinkage estimators, minimaxity is not automatically
achieved but needs to be checked separately.

For the heteroscedastic case, there remain challenging issues on how
much observations with different variances should be shrunk relatively
to each other (e.g., Casella \cite{Cas85}, Morris \cite{Mor83}).
For the empirical Bayes approach (Efron and Morris \cite{EfrMor73}), the
coordinates of $X$ are shrunk directly in proportion to their
variances. But the existing estimators are, in general, non-minimax
(i.e., may have a greater risk than the usual estimator $\delta_0$).
On the other hand, Berger \cite{Ber76} proposed minimax
estimators, including
admissible minimax estimators, such that the coordinates of $X$ are
shrunk inversely in proportion to their variances. But the risk
reduction achieved over $\delta_0$ is insubstantial unless all the
observations have similar variances.

To address the foregoing issues, we develop novel minimax estimators
for multivariate normal means under heteroscedasticity. There are two
central ideas in our approach. The first is to develop a class of
minimax estimators by generalizing a geometric argument essentially in
Stein \cite{Ste56} (see also Brandwein and Strawderman
\cite{BraStr90}). For the
homoscedastic case, the argument shows that $\delta_c^{\mathrm{JS}}$ can be
derived as an approximation to the best linear estimator of the form
$(1-\lambda) X$, where $\lambda$ is a scalar. In fact, the optimal
choice of $\lambda$ in minimizing the risk is $p\sigma^2 /E_\theta(\|
X\|^2)$. Replacing $E_\theta(\|X\|^2)$ by $\|X\|^2$ leads to $\delta
_c^{\mathrm{JS}}$ with $c=p$. This derivation is highly informative, even though
it does not yield the optimal value $c=p-2$.

Our class of minimax estimators are of the linear form $(I - \lambda
A)X$, where $A$ is a nonnegative definite, diagonal matrix indicating
the direction of shrinkage and $\lambda$ is a scalar indicating the
magnitude of shrinkage. The matrix $A$ is open to specification,
depending on the variance matrix $D$ but \textit{not} on the data $X$.
For a fixed $A$, the scalar $\lambda$ is determined to achieve
minimaxity, depending on both $D$ and $X$. Berger's \cite{Ber76} minimax
estimator corresponds to the special choice $A=D^{-1}$, thereby leading
to the unusual pattern of shrinkage discussed above.

The second idea of our approach is to choose $A$ by approximately
minimizing the Bayes risk with a normal prior in our class of minimax
estimators. The Bayes risk is used to measure average risk reduction
for $\theta$ in an elliptical region as in Berger \cite{Ber80,Ber82}. It
turns out that the solution of $A$ obtained by our approximation
strategy has an interesting simple form.
In fact, the coordinates of $X$ are automatically segmented into two
groups, based on their Bayes ``importance'' (Berger \cite
{Ber82}), which is of
the same order as the coordinate variances when the specified prior is
homoscedastic.
The coordinates of high Bayes ``importance'' are shrunk inversely in
proportion to their variances, whereas the remaining coordinates are
shrunk in the direction of the Bayes rule.
This shrinkage pattern may appear paradoxical: it may be expected that
the coordinates of high Bayes ``importance'' are to be shrunk in the
direction of the Bayes rule. But that scheme is inherently aimed at
reducing the Bayes risk under the specified prior and, in general,
fails to achieve minimaxity (i.e., it may lead to even a greater risk
than the usual estimator $\delta_0$).

In addition to simplicity and minimaxity, we further show that the
proposed estimator is scale adaptive in reducing the Bayes risk: it
achieves close to the minimum Bayes risk, with the difference no
greater than the sum of the 4 highest Bayes ``importance'' of the
coordinates of $X$, simultaneously over a scale class of normal priors
(including the specified prior). To our knowledge, the proposed
estimator seems to be the first one with such a property in the general
heteroscedastic case. Previously, in the homoscedastic case, $\delta
_{p-2}^{\mathrm{JS}}$ is known to achieve the minimum Bayes risk up to the sum
of 2 (equal-valued) Bayes ``importance'' of the coordinates over the
scale class of homoscedastic normal priors (Efron and Morris \cite{EfrMor73}).

The rest of this article is organized as follows. Section~\ref{sec2} gives a
review of existing estimators. Section~\ref{sec3} develops the new approach and
studies risk properties of the proposed estimator. Section~\ref{sec4} presents a
simulation study. Section~\ref{sec5} provides concluding remarks. All proofs are
collected in the \hyperref[app]{Appendix}.

\section{Existing estimators}\label{sec2}

We describe a number of existing shrinkage estimators.
See Lehmann and Casella  \cite{LehCas98} for a textbook account and Strawderman
\cite{autokey29} and Morris and Lysy \cite{MorLys12} for recent reviews.
Throughout, $\tr(\cdot)$ denotes the trace and $\lambda_{\max
}(\cdot)$ denotes the largest eigenvalue. Then $\tr(D) = \sum
_{j=1}^p d_j$ and $\lambda_{\max}(D) = \max(d_1,\ldots,d_p)$.

For a Bayes approach, assume the prior distribution: $\theta\sim\N
(0, \gamma I)$,
where $\gamma$ is the prior variance. The Bayes rule is given
componentwise by
$\delta^{\mathrm{Bayes}}_j=\{1- d_j/(d_j + \gamma)\} X_j $.
Then the greater $d_j$ is, the more $X_j$ is shrunk whether $\gamma$
is fixed or estimated from the data.
For the empirical Bayes approach of Efron and Morris \cite{EfrMor73}, $\gamma$
is estimated by the maximum likelihood estimator $\hat\gamma$ such that
%
\begin{eqnarray}\label{EB-iter}
\hat\gamma= \sum_{j=1}^p
\frac{X_j^2- d_j}{(d_j+\hat\gamma)^2} \biggl/ \sum_{j=1}^p
\frac{1}{(d_j+\hat\gamma)^2} .
\end{eqnarray}
Morris \cite{Mor83} suggested the modified estimator
%
\begin{eqnarray}\label{EB}
\delta^{\mathrm{EB}}_j= \biggl( 1- \frac{p-2}{p}
\frac{d_j}{d_j+\hat\gamma
_+} \biggr) X_j .
\end{eqnarray}
In our implementation, the right-hand side of (\ref{EB-iter}) is
computed to update $\hat\gamma$ from the initial guess, $p^{-1} \{
\sum_{j=1}^p (X_j^2-d_j)\}_+$, for up to 100 iterations until the
successive absolute difference in $\hat\gamma$ is $\le$$10^{-4}$, or
$\hat\gamma$ is set to $\infty$ so that $\delta^{\mathrm{EB}}=X$ otherwise.

Alternatively, Xie \textit{et al.} \cite{XieKouBro12} proposed empirical Bayes-type
estimators based on minimizing Stein's \cite{Ste81} unbiased
risk estimate
(SURE) under heteroscedasticity. Their basic estimator is defined
componentwise by
%
\begin{eqnarray}
\delta^{\mathrm{XKB}}_j = \biggl( 1- \frac{d_j}{d_j +\tilde\gamma
} \biggr)
X_j, \label{XKB}
\end{eqnarray}
where $\tilde\gamma$ is obtained by minimizing the SURE of $\delta
^{\mathrm{Bayes}}$, that is, $\operatorname{SURE}(\gamma) =X^\T D \{D+\gamma
I\}^{-1} X
+ 2 \gamma\tr\{D(D+\gamma I)^{-1}\} - \tr(D)$.
In general, the two types of empirical Bayes estimators, $\delta^{\mathrm{EB}}$
and $\delta^{\mathrm{XKB}}$, are non-minimax, as shown in Section~\ref{sec4}.

For a direct extension of $\delta_c^{\mathrm{JS}}$, consider the estimator
$\delta_c^{\mathrm{S}} = (1-c/\|X \|^2) X$ and, more generally,
$\delta_r^{\mathrm{S}} =\{1-r( \|X\|^2 )/\|X \|^2 \} X$, where $c$ is a scalar
constant and $r(\cdot)$ a scalar function. See Lehmann and Casella \cite{LehCas98},
Theorem~5.7, although there are some typos. Both $\delta_c^{\mathrm{S}}$ and
$\delta_r^{\mathrm{S}}$ are spherically symmetric.
The estimator $\delta_c^{\mathrm{S}}$ is minimax provided
%
\begin{eqnarray}\label{S-cond}
0 \le c \le2 \bigl\{\tr(D) - 2 \lambda_{\max} (D) \bigr\},
\end{eqnarray}
and $\delta_r^{\mathrm{S}}$ is minimax provided
$0 \le r(\cdot)\le2 \{\tr(D) - 2 \lambda_{\max} (D) \} \mbox{ and
} r(\cdot) \mbox{ is nondecreasing}$.
No such $c\neq0$ exists unless $ \tr(D) >2 \lambda_{\max} (D) $,
which restricts how much $(d_1, \ldots, d_p)$ can differ from each
other. For example, condition (\ref{S-cond}) fails when $p=10$ and
%
\begin{eqnarray}\label{example}
d_1=40,\qquad  d_2=20,\qquad  d_3=10, \qquad d_4=
\cdots=d_{10} =1,
\end{eqnarray}
because $\tr(D) = 77$ and $\lambda_{\max}(D)=40$.

Berger \cite{Ber76} proposed estimators of the form
$\delta_c^{\mathrm{B}} =\{ I- c D^{-1}/(X^\T D^{-2} X) \} X$ and
$\delta_r^{\mathrm{B}} = \{I- r( X^\T D^{-2} X )/(X^\T D^{-2} X) D^{-1} \} X$,
where $c$ is a scalar constant and $r(\cdot)$ a scalar function.
Then $\delta_c^{\mathrm{B}}$ is minimax provided $0 \le c \le2(p-2)$, and
$\delta_r^{\mathrm{B}}$ is minimax provided
$0 \le r(\cdot) \le2(p-2) \mbox{ and } r(\cdot) \mbox{ is nondecreasing}$,
regardless of differences between $(d_1, \ldots, d_p)$. However, a
striking feature of $\delta_c^{\mathrm{B}}$ and $\delta_r^{\mathrm{B}}$, compared with
$\delta^{\mathrm{EB}}$ and $\delta^{\mathrm{XKB}}$, is that the smaller $d_j$
is, the more $X_j$ is shrunk. For example (\ref{example}), under
$\delta_c^{\mathrm{B}}$, the coordinates $(X_1, X_2, X_3)$ are shrunk only
slightly, whereas $(X_4, \ldots, X_{10})$ are shrunk as if they were
shrunk as a 7-dimensional vector under $\delta_c^{\mathrm{JS}}$. The associated
risk reduction is insubstantial, because the risk of estimating
$(\theta_4,\ldots,\theta_{10})$ is a small fraction of the overall
risk of estimating $\theta$.

Define the positive-part version of $\delta_c^{\mathrm{B}}$ componentwise as
%
\begin{eqnarray}\label{B+}
\bigl(\delta_c^{\mathrm{B}+}\bigr)_j = \biggl( 1-
\frac{c d_j^{-1}}{X^\T D^{-2} X} \biggr)_+ X_j.
\end{eqnarray}
The estimator $\delta_c^{\mathrm{B}+}$ dominates $\delta_c^{\mathrm{B}}$ by Baranchik
\cite{Bar64}, Section~2.5. Berger \cite{Ber85}, Equation (5.32),
stated a different
positive-part estimator, $\delta_r^{\mathrm{B}}$ with $r(t)=\min(p-2, t)$, but
the $j$th component may not be of the same sign as $X_j$.

Given a prior $\theta\sim\N(0,\Gamma)$, Berger \cite
{Ber82} suggested an
approximation of Berger's \cite{Ber80} robust generalized Bayes
estimator as
%
\begin{eqnarray}\label{RB}
\delta^{\mathrm{RB}} = \biggl[I- \min\biggl\{ 1, \frac{p-2}{X^\T(D+\Gamma
)^{-1} X} \biggr\}
D(D+\Gamma)^{-1} \biggr] X.
\end{eqnarray}
The estimator is expected to provide significant risk reduction over
$\delta_0=X$ if the prior is correct and be robust to misspecification
of the prior, but it is, in general, non-minimax.
In the case of $\Gamma=0$, $\delta^{\mathrm{RB}}$ becomes $\{1 - (p-2)/(X^\T
D^{-1} X)\}_+ X$, in the form of spherically symmetric estimators
$\delta^{\mathrm{SS}}_r = \{1 - r(X^\T D^{-1} X)/(X^T D^{-1} X)\}X$, where
$r(\cdot)$ is a scalar function (Bock \cite{Boc75},
Brown \cite{Bro75}). The estimator
$\delta_r^{\mathrm{SS}}$ is minimax provided $0 \le r(\cdot) \le2 \{\tr
(D)/\lambda_{\max}(D)-2\}$ and $r(\cdot)$ is nondecreasing.
Moreover, if $\tr(D) \le 2 \lambda_{\max}(D)$, then $\delta_r^{\mathrm{SS}}$ is non-minimax unless $r(\cdot) =0$.

To overcome the non-minimaxity of $\delta^{\mathrm{RB}}$, Berger
\cite{Ber82}
developed a minimax estimator $\delta^{\mathrm{MB}}$ by combining $\delta
_r^{\mathrm{B}}$, $\delta^{\mathrm{RB}}$, and a minimax estimator of Bhattacharya
\cite{Bha66}.
Suppose that $\Gamma= \diag(\gamma_1, \ldots, \gamma_p)$ and the
indices are sorted such that $d_1^* \ge\cdots\ge d_p^*$, where $d_j^*
= d_j^2/(d_j+\gamma_j)$. Define $\delta^{\mathrm{MB}}$ componentwise as
%
\begin{eqnarray}\label{MB}
\delta^{\mathrm{MB}}_j 
=
X_j - \Biggl[ \frac{1}{d_j^*} \sum_{k=j}^p
\bigl(d_k^*-d_{k+1}^*\bigr) \min\biggl\{1,\frac{ (k-2)_+}{ \sum
_{\ell=1}^k X_\ell^2/(d_\ell+\gamma
_\ell)}
\biggr\} \Biggr]\frac{d_j}{d_j+\gamma_j} X_j,
\end{eqnarray}
%
where $d_{p+1}^*=0$.
In the case of $\Gamma=0$, $\delta^{\mathrm{MB}}$ reduces to the original
estimator of Bhattacharya \cite{Bha66}.
The factor $(k-2)_+$ is replaced by $2(k-2)_+$ in Berger's \cite{Ber82}
original definition of $\delta^{\mathrm{MB}}$,
corresponding to replacing $p-2$ by $2(p-2)$ in $\delta^{\mathrm{RB}}$.
In our simulations, the two versions of $\delta^{\mathrm{MB}}$ somehow yield
rather different risk curves, and so do the corresponding versions of
other estimators.
But there has been limited theory supporting one version over the other.
Therefore, we focus on comparisons of only the corresponding versions
of $\delta^{\mathrm{MB}}$ and other estimators.

\section{Proposed approach}\label{sec3}

We develop a useful approach for shrinkage estimation under
heteroscedasticity, by making explicit how different coordinates are
shrunk differently. The approach not only sheds new light on existing
results, but also lead to new minimax estimators.

\subsection{A sketch}\label{sec3.1}

Assume that $\Sigma=D$ (diagonal) and $Q=I$. Consider estimators of
the linear form
%
\begin{eqnarray}\label{delta-form}
\delta= (I- \lambda A) X = X - \lambda A X,
\end{eqnarray}
where $A$ is a nonnegative definite, diagonal matrix indicating the
\textit{direction} of shrinkage and $\lambda$ is a scalar indicating
the \textit{magnitude} of shrinkage. Both $A$ and $\lambda$ are to be
determined. A sketch of our approach is as follows.
\begin{enumerate}[(iii)]
\item[(i)] For a fixed $A$, the optimal choice of $\lambda$ in
minimizing the risk is
\[
\lambda_{\mathrm{opt}} = \frac{\tr(DA)}{E_\theta(X^\T A^\T A X)}.
\]

\item[(ii)] For a fixed $A$ and a scalar constant $c \ge0$, consider
the estimator
\[
\delta_{A,c} = X - \frac{c}{X^\T A^\T A X} A X .
\]
By Theorem~\ref{th1}, an upper bound on the risk function of $\delta_{A,c}$ is
%
\begin{eqnarray}\label{upper-bound}
R(\delta_{A,c}, \theta) \le\tr(D) + E_\theta\biggl[
\frac{c\{
c-2c^*(D,A)\}}{X^\T A^\T A X} \biggr],
\end{eqnarray}
where $c^*(D,A)=\tr(DA)-2\lambda_{\max}(DA)$.
Requiring the second term to be no greater than 0 shows that if
$c^*(D,A) \ge0$, then $\delta_{A,c}$ is minimax provided
%
\begin{eqnarray}\label{Tan-cond}
0 \le c \le2 c^*(D,A).
\end{eqnarray}
If $c^*(D,A) \ge0$, then the upper bound (\ref{upper-bound}) has a
minimum at $c=c^*(D,A)$.

\item[(iii)] By taking $c=c^*(D,A)$ in $\delta_{A,c}$, consider the estimator
\[
\delta_A = X - \frac{c^*(D,A)}{X^\T A^\T A X} A X
\]
subject to $c^*(D,A)\ge0$, so that $\delta_A$ is minimax by step
(ii). A positive-part estimator dominating $\delta_A$ is defined
componentwise by
%
\begin{eqnarray}\label{A+}
\bigl(\delta_A^+\bigr)_j = \biggl\{1 -
\frac{c^*(D,A) a_j}{X^\T A^\T A X} \biggr\}_+ X_j,
\end{eqnarray}
where $(a_1,\ldots,a_p)$ are the diagonal elements of $A$.
The upper bound (\ref{upper-bound}) on the risk functions of $\delta
_A$ and $\delta^+_{A}$, subject to $c^*(D,A)\ge0$, gives
%
\begin{eqnarray}\label{point-bound}
R(\delta_A, \theta) \le\tr(D) - E_\theta\biggl\{
\frac
{{c^*}^2(D,A)}{X^\T A^\T A X} \biggr\}.
\end{eqnarray}
We propose to choose $A$ based on some optimality criterion, such as
minimizing the Bayes risk with a normal prior centered at 0 (Berger \cite{Ber82}).
\end{enumerate}
Further discussions of steps (i)--(iii) are provided in Sections~\ref{sec3.2}--\ref{sec3.3}.

\subsection{Constructing estimators: Steps (i)--(ii)}\label{sec3.2}

We first develop steps (i)--(ii) for the general problem where neither
$\Sigma$ nor $Q$ may be diagonal. The results can be as concisely
stated as those just presented for the canonical problem where $\Sigma
$ is diagonal and $Q=I$. Such a unification adds to the attractiveness
of the proposed approach.

Consider estimators of the form (\ref{delta-form}), where $A$ is not
necessarily diagonal, but
%
\begin{eqnarray}\label{A-cond}
A \Sigma\mbox{ is nonnegative definite.}
\end{eqnarray}
Condition (\ref{A-cond}) is invariant under a linear transformation.
To see this, let $B$ be a nonsingular matrix and $\Sigma^*=B \Sigma
B^\T$ and $A^*=B A B^{-1}$. For the transformed problem of estimating
$\theta^*=B \theta$ based on $X^*=B X$ with variance matrix $\Sigma
^*$, the transformed estimator from (\ref{delta-form}) is $\delta^* =
X^* - \lambda A^* X^*$. The application of condition (\ref{A-cond}) to
$\delta^*$ says that $A^* \Sigma^* = B A \Sigma B^\T$ is nonnegative
definite and therefore is equivalent to (\ref{A-cond}) itself. For the
canonical problem where $\Sigma=D$ (diagonal), condition (\ref
{A-cond}) only requires that $AD$ is nonnegative definite, allowing $A$
to be non-diagonal. On the other hand, it seems intuitively appropriate
to restrict $A$ to be diagonal. Then condition (\ref{A-cond}) is
equivalent to saying that $A$ is nonnegative definite (and diagonal),
which is the condition introduced on $A$ in the sketch in Section~\ref{sec3.1}.

The risk of an estimator of the form (\ref{delta-form}) is
\begin{eqnarray*}
&& E_\theta\bigl\{ (X-\theta-\lambda A X)^\T Q (X-\theta-
\lambda A X) \bigr\}
\\
&&\quad =  E_\theta\bigl\{ (X-\theta)^\T Q (X-\theta) \bigr\} +
\lambda^2 E_\theta\bigl(X^\T A^\T Q A
X \bigr) - 2 \lambda E_\theta\bigl\{ (X-\theta)^\T Q A X
\bigr\} .
\end{eqnarray*}
For a fixed $A$, the optimal $\lambda$ in minimizing the risk is
\begin{eqnarray*}
\lambda_{\mathrm{opt}} = \frac{E_\theta\{ (X-\theta)^\T Q A X \}
}{E_\theta(X^\T A^\T Q A X )} = \frac{\tr(\Sigma Q A)}{E_\theta
(X^\T A^\T Q A X )}.
\end{eqnarray*}
Replacing $E_\theta(X^\T A^\T Q A X )$ by $X^\T A^\T Q A X$ and $\tr
(\Sigma QA)$ by a scalar constant $c\ge0$ leads to the estimator
\[
\delta_{A,c} = X - \frac{c}{X^\T A^\T Q A X} A X .
\]
For a generalization, replacing $c$ by $r(X^\T A^\T Q A X)$ with a
scalar function $r(\cdot)\ge0$ leads to the estimator
\[
\delta_{A,r} = X - \frac{r(X^\T A^\T Q A X)}{X^\T A^\T Q A X} A X .
\]
We provide in Theorem~\ref{th1} an upper bound on the risk function of $\delta_{A,r}$.
%
\begin{Theorem}\label{th1}
Assume that $r(\cdot)$ almost differentiable
(Stein \cite{Ste81}). If (\ref{A-cond}) holds and
$r(\cdot) \ge0$ is
nondecreasing, then for each $\theta$,
%
\begin{eqnarray}\label{r-upper-bound}
R(\delta_{A,r}, \theta) \le\tr(\Sigma Q) + E_\theta\biggl[
\frac
{r\{r-2c^*(\Sigma,Q,A)\}}{X^\T A^\T Q A X} \biggr],
\end{eqnarray}
where $r=r(X^\T A^\T Q A X)$ and $c^*(\Sigma,Q,A)=\tr(A \Sigma
Q)-\lambda_{\max}(A \Sigma Q + \Sigma A^\T Q )$. Taking $r(\cdot
)\equiv c \ge0$ in (\ref{r-upper-bound}) gives an upper bound on
$R(\delta_{A,c}, \theta)$.
\end{Theorem}

Requiring the second term in the risk upper bound (\ref
{r-upper-bound}) to be no greater than 0 leads to a sufficient
condition for $\delta_{A,r}$ to be minimax.
%
\begin{Corollary}\label{cor1} If (\ref{A-cond}) holds and $c^*(\Sigma
,Q,A)\ge0$, then $\delta_{A,r}$ is minimax provided
%
\begin{eqnarray}\label{Tan-cond2}
0 \le r(\cdot) \le2c^*(\Sigma,Q,A)\quad  \mbox{and} \quad r(\cdot) \mbox{ is
nondecreasing}.
\end{eqnarray}
Particularly, $\delta_{A,c}$ is minimax provided $0 \le c \le
2c^*(\Sigma,Q,A)$.
\end{Corollary}

For the canonical problem, inequality (\ref{r-upper-bound}) and
condition (\ref{Tan-cond2}) for $\delta_{A,c}$ give respectively
(\ref{upper-bound}) and (\ref{Tan-cond}). These results generalize
the corresponding ones for $\delta_c^{\mathrm{S}}$ and $\delta_c^{\mathrm{B}}$ in
Section~\ref{sec2}, by the specific choices $A=I$ or $D^{-1}$. The
generalization also
holds if $c$ is replaced by a scalar function $r(\cdot)>0$. In fact,
condition (\ref{Tan-cond2}) reduces to Baranchik's \cite{Bar70}
condition in
the homoscedastic case.

If $c^*(\Sigma,Q,A)\ge0$, then the risk upper bound (\ref
{r-upper-bound}) has a minimum at $r(\cdot) \equiv c = c^*(\Sigma,
Q,A)$. As a result, consider the estimator
\begin{eqnarray*}
\delta_A = X - \frac{c^*(\Sigma,Q,A)}{X^\T A^\T Q A X} A X,
\end{eqnarray*}
which is minimax provided $c^*(\Sigma,Q,A)\ge0$. If $A=Q^{-1}\Sigma
^{-1}$ (Berger \cite{Ber76}), then $c^*(\Sigma,Q,A)=p-2$
and, by the proof of
Theorem~\ref{th1} in the \hyperref[app]{Appendix}, the risk upper bound (\ref{r-upper-bound})
becomes exact for $\delta_{A,c}$. Therefore, for $A=Q^{-1}\Sigma
^{-1}$, the estimator $\delta_A=\delta_{A,p-2}$ is uniformly best in
the class $\delta_{A,c}$, in agreement with the result that $\delta
_{p-2}^{\mathrm{JS}}$ is uniformly best among $\delta_c^{\mathrm{JS}}$ in the
homoscedastic case.

The estimator $\delta_A$ has desirable properties of invariance.
First, $\delta_A$ is easily shown to be invariant under a
multiplicative transformation $A \mapsto aA$ for a scalar $a >0$.
Second, $\delta_A$ is invariant under a linear transformation of the
inference problem. Similarly as discussed below (\ref{A-cond}), let
$B$ be a nonsingular matrix and $\Sigma^*=B \Sigma B^\T$, $Q^*={B^\T
}^{-1} Q B^{-1}$, and $A^*=B A B^{-1}$.
For the transformed problem of estimating $\theta^*=B \theta$ based
on $X^*=B X$, the transformed estimator from $\delta_A$ is $ X^* - \{
c^*(\Sigma,Q,A)/({X^*}^\T{A^*}^\T Q^* A^* X^*)\} A^* X^*$, whereas
the application of $\delta_A$ is
$ X^* - \{c^*(\Sigma^*,Q^*,A^*)/({X^*}^\T{A^*}^\T Q^* A^* X^*)\} A^*
X^*$. The two estimators are identical because $A^*\Sigma^*Q^* = B A
\Sigma Q B^{-1}$, $\Sigma^* {A^*}^\T Q^* =B \Sigma A^\T Q B^{-1}$, and
hence $c^*(\Sigma^*,Q^*,A^*)=c^*(\Sigma,Q,A)$.

Finally, we present a positive-part estimator dominating $\delta_A$ in
the case where both $A\Sigma$ and $QA$ are symmetric, that is,
%
\begin{eqnarray}\label{A-cond2}
A\Sigma=\Sigma A^\T\quad \mbox{and}\quad  QA =A^\T Q .
\end{eqnarray}
Similarly to (\ref{A-cond}), it is easy to see that this condition is
invariant under a linear transformation.
Condition (\ref{A-cond2}) is trivially true if $\Sigma$, $Q$, and $A$
are diagonal.
In the \hyperref[app]{Appendix}, we show that (\ref{A-cond2}) holds if and only if
there exists a nonsingular matrix $B$ such that $Q=B^\T B$, $\Sigma
=B^{-1}D {B^\T}^{-1}$, and $A=B^{-1} A^* B$, where $D$ and $A^*$ are
diagonal and the diagonal elements of $D$ or $A^*$ are, respectively,
the eigenvalues of $\Sigma Q$ or $A$. In the foregoing notation,
$\Sigma^*=D$ and $Q^*=I$. For the problem of estimating $\theta^*=B
\theta$ based on $X^*=B X$, consider the estimator $\eta= X- \{
c^*(D,A^*)/({X^*}^\T{A^*}^\T A^* X^*) \} A^* X$ and the positive-part
estimator $\eta^+$ with the $j$th component,
\begin{eqnarray*}
\biggl\{ 1- \frac{c^*(D,A^*)}{{X^*}^\T{A^*}^\T A^* X^*} a_j^* \biggr
\}_+ X^*_j
,
\end{eqnarray*}
where $(a_1^*,\ldots,a_p^*)$ are the diagonal elements of $A^*$.
The estimator $\eta^+$ dominates $\eta$ by a simple extension of
Baranchik \cite{Bar64}, Section~2.5. By a transformation back to
the original
problem, $\eta$ yields $\delta_A$, whereas $\eta^+ $ yields
\begin{eqnarray*}
\delta_A^+ = B^{-1} \diag\biggl[ \biggl\{1-
\frac{c^*(\Sigma
,Q,A)}{X^\T A^\T Q AX}a^*_1 \biggr\}_+, \ldots, \biggl\{1-
\frac
{c^*(\Sigma,Q,A)}{X^\T A^\T Q AX} a^*_p \biggr\}_+ \biggr] B X.
\end{eqnarray*}
Then $\delta_A^+$ dominates $\delta_A$. Therefore, (\ref
{r-upper-bound}) also gives an upper bound on the risk of $\delta
_A^+$, with $r(\cdot) \equiv c^*(\Sigma, Q,A)$, even though $\delta
_A^+$ is not of the form $\delta_{A,r}$.

In practice, a matrix $A$ satisfying (\ref{A-cond2}) can be specified
in two steps. First, find a nonsingular matrix $B$ such that $Q=B^\T B$
and $\Sigma=B^{-1}D {B^\T}^{-1}$, where $D$ is diagonal.
Second, pick a diagonal matrix $A^*$ and define $A= B^{-1} A^* B$. The
first step is always feasible by taking $B=OC$, where $C$ is a
nonsingular matrix such that $Q= C^\T C$ and $O$ is an orthogonal
matrix $O$ such that $O (C\Sigma C^\T) O^\T$ is diagonal. Given
$(\Sigma,Q)$ and $D$, it can be shown that $A$ and $\delta_A^+$
depend on the choice of $A^*$, but not on that of $B$, provided that
$a^*_j=a^*_k$ if $d_j=d_k$ for any $j, k=1,\ldots,p$. In the canonical
case where $\Sigma=D$ and $Q=I$, this condition amounts to saying that
any coordinates of $X$ with the same variances should be shrunk in the
same way.

\subsection{Constructing estimators: Step (iii)}\label{sec3.3}

Different choices of $A$ lead to different estimators $\delta_A$ and
$\delta_A^+$. We study how to choose $A$, depending on $(\Sigma, Q)$
but \textit{not} on $X$, to approximately optimize risk reduction
while preserving minimaxity for $\delta_A$. The estimator $\delta
_A^+$ provides even greater risk reduction than $\delta_A$. We focus
on the canonical problem where $\Sigma=D$ (diagonal) and $Q=I$.
Further, we restrict $A$ to be diagonal and nonnegative definite.

As discussed in Berger \cite{Ber80}, any estimator can
have significantly
smaller risk than $\delta_0=X$ only for $\theta$ in a specific
region. Berger \cite{Ber80,Ber82} considered the
situation where significant
risk reduction is desired for an elliptical region
%
\begin{eqnarray}\label{region}
\bigl\{\theta\dvt  (\theta-\mu)^\T\Gamma^{-1} (\theta-\mu) \le
p\bigr\},
\end{eqnarray}
with $\mu$ and $\Gamma$ the prior mean and prior variance matrix. See
$\delta^{\mathrm{RB}}$ and $\delta^{\mathrm{MB}}$ reviewed in Section~\ref{sec2}. To measure
average risk reduction for $\theta$ in region (\ref{region}),
Berger \cite{Ber82} used the Bayes risk with the normal
prior $\theta\sim\N(\mu
,\Gamma)$. For simplicity, assume throughout that $\mu=0$ and $\Gamma
=\diag(\gamma_1,\ldots,\gamma_p)$ is diagonal.

We adopt Berger's \cite{Ber82} ideas of specifying an elliptical
region and
using the Bayes risk to quantify average risk reduction in this region.
We aim to find $A$, subject to $c^*(D,A)\ge0$, minimizing the Bayes
risk of $\delta_A$ with the prior $\pi_\Gamma$, $\theta\sim\N
(0,\Gamma)$,
\begin{eqnarray*}
R(\delta_A, \pi_\Gamma) = E^{\pi_\Gamma} E_\theta
\bigl( \|\delta_A-\theta\|^2 \bigr),
\end{eqnarray*}
where $E^{\pi_\Gamma}$ denotes the expectation with respect to the
prior $\pi_\Gamma$. Given $A$, the risk $R(\delta_A, \pi_\Gamma)$
can be numerically evaluated. A simple Monte Carlo method is to
repeatedly draw $\theta\sim\N(0,\Gamma)$ and $X |\theta\sim\N
(\theta, D)$ and then take the average of $\| \delta_A(X) - \theta\|
^2$. But it seems difficult to literally implement the foregoing
optimization. Alternatively, we develop a simple method for choosing
$A$ by two approximations.

First, if $c^*(D,A)\ge0$, then taking the expectation of both sides of
(\ref{point-bound}) with respect to the prior $\pi_\Gamma$ gives an
upper bound on the Bayes risk of $\delta_A$:
%
\begin{eqnarray}\label{bayes-bound}
R(\delta_A, \pi_\Gamma) \le\tr(D) - E^m \biggl
\{ \frac
{{c^*}^2(D,A)}{X^\T A^\T A X} \biggr\},
\end{eqnarray}
where $E^m$ denotes the expectation with respect to the marginal
distribution of $X$ in the Bayes model, that is, $X \sim\N(0,D+\Gamma
)$. An approximation strategy for choosing $A$ is to minimize the upper
bound (\ref{bayes-bound}) on the Bayes risk or to maximize the second term.
The expectation $E^m\{(X^\T A^\T A X)^{-1}\}$ can be evaluated as a
1-dimensional integral by results on inverse moments of quadratic forms
in normal variables (e.g., Jones \cite{Jon86}). But the
required optimization
problem remains difficult.

Second, approximations can be made to the distribution of the quadratic
form $X^\T A^\T A X$. Suppose that $X^\T A^\T A X$ is approximated with
the same mean by $\{\sum_{j=1}^p (d_j+\gamma_j)a_j^2 \} \chi^2_p/p$,
where $\chi^2_p$ is a chi-squared variable with $p$ degrees of
freedom. Then $E^m\{(X^\T A^\T A X)^{-1}\}$ is approximated by $\{
p/(p-2)\} \{\sum_{j=1}^p (d_j+\gamma_j)a_j^2\}^{-1}$. We show in the
\hyperref[app]{Appendix} that this approximation gives a valid lower bound:
%
\begin{eqnarray}\label{bayes-bound2}
E^m \biggl( \frac{1}{X^\T A^\T A X} \biggr) \ge\frac{p}{p-2} \cdot
\frac{1}{ \sum_{j=1}^p (d_j+\gamma_j)a_j^2} .
\end{eqnarray}
A direct application of Jensen's inequality shows that $E^m\{(X^\T A^\T
A X)^{-1}\} \ge\{\sum_{j=1}^p (d_j+\gamma_j)a_j^2\}^{-1}$.
But the lower bound (\ref{bayes-bound2}) is strictly tighter and
becomes exact when $(d_1+\gamma_1)a_1^2 = \cdots=(d_p+\gamma
_p)a_p^2$. No simple bounds such as (\ref{bayes-bound2}) seem to hold
if more complicated approximations (e.g., Satterthwaite \cite{SAT46}) are used.

Combining (\ref{bayes-bound}) and (\ref{bayes-bound2}) shows that if
$c^*(D,A)\ge0$, then
%
\begin{eqnarray}\label{bayes-bound3}
R(\delta_A, \pi_\Gamma) \le\tr(D) - \frac{p}{p-2}
\cdot\frac
{{c^*}^2(D,A)}{ \sum_{j=1}^p (d_j+\gamma_j)a_j^2} .
\end{eqnarray}
Notice that $\delta_A$ is invariant under a multiplicative
transformation $A \mapsto a A$ for a scalar $a >0$, and so is the upper
bound (\ref{bayes-bound3}).
Our strategy for choosing $A$ is to minimize the upper bound (\ref
{bayes-bound3}) subject to $c^*(D,A)\ge0$ or, equivalently, to solve
the constrained optimization problem:
%
\begin{eqnarray}\label{opt}
&&\max_A \quad  c^*(D,A)=\sum_{j=1}^p
d_j a_j - 2 \max_{j=1,\ldots
,p}
d_j a_j \nonumber
\\[-8pt]\\[-8pt]
&&\quad \mbox{subject to} \quad  \sum_{j=1}^p
(d_j+ \gamma_j) a_j^2 = \mbox{fixed}.
\nonumber
\end{eqnarray}
The condition $c^*(D,A)\ge0$ is dropped, because
for $p\ge3$, the achieved maximum is at least $c^*(D,
aD^{-1})=a(p-2)>0$ for some scalar $a>0$. In spite of the
approximations used in our approach, Theorem~\ref{th2} shows that not only the
problem (\ref{opt}) admits a non-iterative solution, but also the
solution has a very interesting interpretation.
For convenience, assume thereafter that the indices are sorted such that
$d_1^2/(d_1+\gamma_1) \ge d_2^2/(d_2+\gamma_2) \ge\cdots\ge
d_p^2/(d_p+\gamma_p)$.
%
\begin{Theorem}\label{th2} Assume that $p\ge3$, $D=\diag(d_1, \ldots,
d_p)$ with $d_j >0$ and $\Gamma= \diag(\gamma_1, \ldots, \gamma
_p)$ with $\gamma_j\ge0$ ($j=1,\ldots,p$). For problem (\ref{opt}),
assume that $A=\diag(a_1, \ldots, a_p)$ with $a_j \ge0$ ($j=1,\ldots
,p$) and $\sum_{j=1}^p (d_j+\gamma_j) a_j^2= \sum_{j=1}^p
d_j^2/(d_j+\gamma_j)$, satisfied by $a_j=d_j/(d_j+\gamma_j)$.
Then the following results hold.
\begin{enumerate}[(iii)]
\item[(i)] There exists a \textit{unique} solution, $A^\dag= \diag
(a_1^\dag, \ldots, a_p^\dag)$, to problem (\ref{opt}).

\item[(ii)] Let $\nu$ be the largest index such that $d_\nu a^\dag
_\nu= \max(d_1a_1^\dag,\ldots,d_pa_p^\dag)$. Then $\nu\ge3$,
$d_1a_1^\dag= \cdots= d_\nu a_\nu^\dag> d_ja_j^\dag$ for $j \ge
\nu+1$, and
\begin{eqnarray*}
a_j^\dag& =& K_\nu\Biggl(\sum
_{k=1}^\nu\frac{d_k+\gamma
_k}{d_k^2} \Biggr)^{-1}
\frac{\nu-2}{d_j}\qquad  (j=1,\ldots, \nu),
\\
a_j^\dag& =& K_\nu\frac{d_j}{d_j+\gamma_j}\qquad  (j=\nu+1,
\ldots, p),
\end{eqnarray*}
where $K_\nu= \{\sum_{j=1}^p d_j^2/(d_j+\gamma_j)\}^{1/2} M_\nu
^{-1/2}$ and
\begin{eqnarray*}
M_\nu= \frac{(\nu-2)^2}{ \sum_{j=1}^\nu\vfrac{d_j+\gamma
_j}{d_j^2}} + \sum_{j=\nu+1}^p
\frac{d_j^2}{d_j+\gamma_j} .
\end{eqnarray*}
The achieved maximum value, $c^*(D,A^\dag)$, 
is $K_\nu M_\nu\ (>0)$.
\item[(iii)] The resulting estimator $\delta_{A^\dag}$ is minimax.
\end{enumerate}
\end{Theorem}

We emphasize that, although $A$ can be considered a tuning parameter,
the solution $A^\dag$ is \textit{data independent}, so that $\delta
_{A^\dag}$ is automatically minimax. If a data-dependent choice of $A$
were used, minimaxity would not necessarily hold. This result is
achieved both because each estimator $\delta_A$ with $c^*(D,A)\ge0$
is minimax and because a global criterion (such as the Bayes risk) is
used, instead of a pointwise criterion (such as the frequentist risk at
the unknown $\theta$), to select $A$. By these considerations, our
approach differs from the usual exercise of selecting a tuning
parameter in a data-dependent manner for a class of candidate estimators.

There is a remarkable property of monotonicity for the sequence $(M_3,
M_4, \ldots, M_p)$, which underlies the uniqueness of $\nu$ and
$A^\dag$.
%
\begin{Corollary}\label{cor2} The sequence $(M_3, M_4, \ldots, M_p)$ is
nonincreasing: for $3 \le k \le p-1$, $M_k \ge M_{k+1}$, where the
equality holds if and only if
\[
\frac{k-2}{  \sum_{j=1}^k \vfrac{d_j+\gamma_j}{d_j^2} } = \frac
{d_{k+1}^2}{d_{k+1}+\gamma_{k+1}} .
\]
The condition $d_\nu a_\nu^\dag> d_{\nu+1} a_{\nu+1}^\dag$ is
equivalent to saying that the left side is greater than the right-hand
side in the above expression for $k=\nu$. Therefore, $\nu$ is the
smallest index $3 \le k \le p-1$ with this property, and $M_\nu>
M_{\nu+1}$.
\end{Corollary}

The estimator $\delta_{A^\dag}$ is invariant under scale
transformations of $A^\dag$. Therefore, the constant $K_\nu$ can be
dropped from the expression of $A^\dag$ in Theorem~\ref{th1}.
%
\begin{Corollary}\label{cor3} The solution $A^\dag=\diag(a_1^\dag, \ldots,
a_p^\dag)$ can be rescaled such that
%
\begin{eqnarray}
\label{sol1}a_j^\dag& =& \Biggl(\sum_{k=1}^\nu
\frac{d_k+\gamma_k}{d_k^2} \Biggr)^{-1} \frac{\nu-2}{d_j}\qquad
(j=1,\ldots,\nu),
\\
\label{sol2}a_j^\dag& =& \frac{d_j}{d_j + \gamma_j} \qquad (j=\nu+1, \ldots, p).
\end{eqnarray}
Then $c^*(D, A^\dag) = \sum_{j=1}^p {a_j^{\dag}}^2 (d_j+\gamma_j) =
M_\nu$. Moreover, it holds that
%
\begin{eqnarray}\label{sol-ineq}
a_j^\dag\le\frac{d_j}{d_j + \gamma_j} \qquad (j=1,\ldots,\nu).
\end{eqnarray}
The estimator $\delta_{A^\dag}$ can be expressed as
%
\begin{eqnarray}\label{delta-A}
\delta_{A^\dag} 
= X -
\frac{  \sum_{j=1}^p {a_j^{\dag}}^2 (d_j+\gamma_j) }{  \sum_{j=1}^p
{a_j^{\dag}}^2 X_j^2} A^\dag X.
\end{eqnarray}
\end{Corollary}

The foregoing results lead to a simple algorithm for solving problem
(\ref{opt}):
\begin{enumerate}[(iii)]
\item[(i)] Sort the indices such that $d_1^2/(d_1+\gamma_1) \ge
\cdots\ge d_p^2/(d_p+\gamma_p)$.

\item[(ii)] Take $\nu$ to be the smallest index $k$ (corresponding to
the largest $M_k$) such that $3 \le k \le p-1$ and
\[
\frac{k-2}{  \sum_{j=1}^k \vfrac{d_j+\gamma_j}{d_j^2} } > \frac
{d_{k+1}^2}{d_{k+1}+\gamma_{k+1}} ,
\]
or take $\nu=p$ if there exists no such $k$.

\item[(iii)] Compute $(a_1^\dag, \ldots, a_p^\dag)$ by (\ref
{sol1})--(\ref{sol2}).
\end{enumerate}
This algorithm is guaranteed to find the (unique) solution to problem
(\ref{opt}) by a fixed number of numerical operations. No iteration or
convergence diagnosis is required. Therefore, the algorithm is exact
and non-iterative, in contrast with usual iterative algorithms for
nonlinear, constrained optimization.

The estimator $\delta_{A^\dag}$ has an interesting interpretation. By
(\ref{sol1})--(\ref{sol2}), there is a dichotomous segmentation in
the shrinkage direction of the coordinates of $X$ based on
$d_j^*=d_j^2/(d_j+\gamma_j)$. This quantity $d_j^*$ is said to reflect
the Bayes ``importance'' of $\theta_j$, that is, the amount of
reduction in Bayes risk obtainable in estimating $\theta_j$ in
Berger \cite{Ber82}. The coordinates with high $d_j^*$
are shrunk inversely in
proportion to their variances $d_j$ as in Berger's \cite{Ber76} estimator
$\delta_c^{\mathrm{B}}$, whereas the coordinates with low $d_j^*$ are shrunk in
the direction of the Bayes rule. Therefore,
$\delta_{A^\dag}$ mimics the Bayes rule to reduce the Bayes risk,
except that $\delta_{A^\dag}$ mimics $\delta_c^{\mathrm{B}}$ for some
coordinates of highest Bayes ``importance'' in order to achieve
minimaxity. In fact, by inequality (\ref{sol-ineq}), the relative
shrinkage, $a_j^\dag/\{d_j/(d_j+\gamma_j)\}$, of each $X_j$
($j=1,\ldots,\nu$) in $\delta_{A^\dag}$ versus the Bayes rule is
always no greater than that of $X_k$ ($k=\nu+1,\ldots,p$).

The expression (\ref{delta-A}) suggests that there is a close
relationship in beyond the shrinkage direction between $\delta_{A^\dag
}$ and the Bayes rule under the Bayes model, $X \sim\N(0, D+\Gamma
)$. In this case, $E^m( \sum_{j=1}^p {a_j^{\dag}}^2 X_j^2 ) = \sum
_{j=1}^p {a_j^{\dag}}^2 (d_j +\gamma_j)$, and hence $\delta_{A^\dag
}$ behaves similarly to $X - A^\dag X$. Therefore, \textit{on average}
under the Bayes model, the coordinates of $X$ are shrunk in $\delta
_{A^\dag}$ the same as in the Bayes rule, except that some coordinates
of highest Bayes ``importance'' are shrunk no greater than in the Bayes rule.
While this discussion seems heuristic, we provide in Section~\ref{sec3.4} a
rigorous analysis of the Bayes risk of $\delta_{A^\dag}$, compared
with that of the Bayes rule.

We now examine $\delta_{A^\dag}$ for two types of priors: $\gamma
_1=\cdots=\gamma_p = \gamma$ and $\gamma_j = \gamma d_j $ ($j=1,
\ldots, p$), referred to as the homoscedastic and heteroscedastic priors.
For both types, $(d_1^*,\ldots,d_p^*)$ are of the same order as the
variances $(d_1,\ldots,d_p)$.
Recall that $\delta_A$ is invariant under a multiplicative
transformation of $A$.
For both the homoscedastic prior with $\gamma=0$ and the
heteroscedastic prior \textit{regardless} of $\gamma\ge0$,
the solution $A^\dag=\diag(a_1^\dag, \ldots, a_p^\dag)$ can be
rescaled such that
\begin{eqnarray*}
a_j^\dag&=& \Biggl(\sum_{k=1}^\nu
d_k^{-1} \Biggr)^{-1} \frac{\nu
-2}{d_j} \qquad (j=1,
\ldots, \nu),
\\
a_j^\dag&=& 1 \qquad (j=\nu+1,\ldots, p).
\end{eqnarray*}
Denote by $A^\dag_0$ this rescaled matrix $A^\dag$, corresponding to
$\Gamma=0$.
Then coordinates with high variances are shrunk inversely in proportion
to their variances, whereas coordinates with low variances are shrunk
symmetrically. For $\Gamma=0$, the proposed method has a purely
frequentist interpretation: it seeks to minimize the upper bound (\ref
{bayes-bound3}) on the pointwise risk of $\delta_A$ at \mbox{$\theta=0$}.

For the homoscedastic prior with $\gamma\to\infty$, the proposed
method is then to minimize the upper bound (\ref{bayes-bound3}) on the
Bayes risk of $\delta_A$ with an extremely flat, homoscedastic prior.
As $\gamma\to\infty$, the solution $A^\dag$ can be rescaled such that
\begin{eqnarray*}
a_j^\dag&=& \Biggl(\sum_{k=1}^\nu
d_k^{-2} \Biggr)^{-1} \frac{\nu
-2}{d_j} \qquad (j=1,
\ldots, \nu),
\\
a_j^\dag&=& d_j\qquad  (j=\nu+1,\ldots, p).
\end{eqnarray*}
Denote by $A^\dag_\infty$ this rescaled matrix $A^\dag$.
Then coordinates with low (or high) variances are shrunk directly (or
inversely) in proportion to their variances.
The direction $A^\dag_\infty$ can also be\vadjust{\goodbreak} obtained by using a fixed
prior in the form $\gamma_j = \gamma d_1 - d_j$ ($j=1,\ldots,p$) for
arbitrary $\gamma\ge1$, where $d_1 = \max_{j=1,\ldots,p} d_j$.

Finally, in the homoscedastic case ($d_1=\cdots=d_p=\sigma^2$), if
the prior is also homoscedastic ($\gamma_1=\cdots=\gamma_p=\gamma
$), then $\nu=p$, $a_1^\dag= \cdots= a_p^\dag$, and $\delta
_{A^\dag}$ reduces to the James--Stein estimator $\delta_{p-2}^{\mathrm{JS}}$,
\textit{regardless} of $\sigma^2$ and $\gamma$.

\subsection{Evaluating estimators}\label{sec3.4}

The estimator $\delta_{A^\dag}$ is constructed by minimizing the
upper bound (\ref{bayes-bound3}) on the Bayes risk subject to
minimaxity. In addition to simplicity, interpretability, and minimaxity
demonstrated for $\delta_{A^\dag}$, it remains important to further
study risk properties of $\delta_{A^\dag}$ and show that $\delta
_{A^\dag}$ can provide effective risk reduction over $\delta_0=X$.
Write $\delta_{A^\dag} = \delta_{A^\dag(\Gamma)}$ whenever needed
to make explicit the dependency of $A^\dag$ on $\Gamma$.

First, we study how close the Bayes risk of $\delta_{A^\dag(\Gamma
)}$ can be to that of the Bayes rule, which is the smallest possible
among \textit{all} estimators including non-minimax ones, under the
prior $\pi_\Gamma$, $\theta\sim\N(0, \Gamma)$. The Bayes rule
$\delta_\Gamma^{\mathrm{Bayes}}$ is given componentwise by $(\delta
^{\mathrm{Bayes}}_\Gamma)_j = \{1- d_j/(d_j+\gamma_j)\} X_j$, with the Bayes risk\vspace*{1pt}
\begin{eqnarray*}
R\bigl(\delta^{\mathrm{Bayes}}_\Gamma, \pi_\Gamma\bigr) = \tr(D) -
\sum_{j=1}^p d^*_j ,
\end{eqnarray*}
where $d^*_j =d_j^2/(d_j+\gamma_j)$, indicating the Bayes
``importance'' of $\theta_j$ (Berger \cite{Ber82}).
The upper bound (\ref{bayes-bound3}) on the Bayes risk of $\delta
_{A^\dag(\Gamma)}$ gives\vspace*{1pt}
%
\begin{eqnarray}\label{bayes-bound4}
R\{\delta_{A^\dag(\Gamma)}, \pi_\Gamma\} \le\tr(D) - \frac
{p}{p-2}
M_\nu= \tr(D) - \frac{p}{p-2} \Biggl\{ \frac{(\nu
-2)^2}{ \sum_{j=1}^\nu{d_j^*}^{-1}} + \sum
_{j=\nu+1}^p d^*_j \Biggr\} ,
\end{eqnarray}
because $c^*(D,A^\dag) = \sum_{j=1}^p (d_j+\gamma_j) {a_j^\dag}^2 =
M_\nu$  and hence ${c^*}^2(D, A^\dag) /  \{\sum
_{j=1}^p (d_j+\gamma_j) {a_j^\dag}^2\} = M_\nu$ by Corollary~\ref{cor3}. It appears that the
difference between $R\{\delta_{A^\dag(\Gamma)}, \pi_\Gamma\}$ and
$R(\delta^{\mathrm{Bayes}}_\Gamma,\allowbreak  \pi_\Gamma)$ tends to be large if $\nu$
is large. But $d_1^* \ge\cdots\ge d_\nu^*$ cannot differ too much
from each other because by Corollary~\ref{cor1},\vspace*{1pt}
\[
k-2 \le\sum_{j=1}^k \frac{d_{k+1}^*}{d_j^*}
\le k \qquad (k=3,\ldots,\nu-1).
\]
Then the difference between $R\{\delta_{A^\dag(\Gamma)}, \pi_\Gamma
\}$ and $R(\delta^{\mathrm{Bayes}}_\Gamma, \pi_\Gamma)$ should be limited
even if $\nu$ is large. A careful analysis using these ideas leads to
the following result.
%
\begin{Theorem}\label{th3} Suppose that the prior is $\theta\sim\N
(0,\Gamma)$. If $\nu=3$, then
%
\begin{eqnarray}
\label{tight-bound1}R\{\delta_{A^\dag(\Gamma)}, \pi_\Gamma\} & \le&\tr(D) - \sum
_{j=3}^p d_j^* + \Biggl(
d_3^* - \frac{2}{p-2} \sum_{j=4}^p
d_j^* -\frac{ p}{p-2} \frac{d_3^*}{3} \Biggr)
\\
\label{loose-bound1}& \le&\tr(D) - \sum_{j=3}^p
d_j^* + \frac{2}{3} d_3^*.
\end{eqnarray}
If $\nu\ge4$, then
%
\begin{eqnarray}
\label
{tight-bound2}R\{\delta_{A^\dag(\Gamma)}, \pi_\Gamma\} & \le&\tr(D) - \sum
_{j=3}^p d_j^* + \Biggl(
d_3^* + d_4^* - \frac{2}{p-2} \sum
_{j=5}^p d_j^* -\frac{ 4 p}{p-2}
\frac{d_\nu^*}{\nu} \Biggr)
\\
\label{loose-bound2}& \le&\tr(D) - \sum_{j=3}^p
d_j^* + \bigl( d_3^* + d_4^* \bigr).
\end{eqnarray}
Throughout, an empty summation is 0.
\end{Theorem}

There are interesting implications of Theorem~3. By (\ref
{loose-bound1}) and (\ref{loose-bound2}),
%
\begin{eqnarray}\label
{bayes-close}
R\{\delta_{A^\dag(\Gamma)}, \pi_\Gamma\} \le R\bigl(\delta
^{\mathrm{Bayes}}_\Gamma, \pi_\Gamma\bigr) + \bigl(d_1^*+d_2^*+d_3^*+d_4^*
\bigr).
\end{eqnarray}
Then $\delta_{A^\dag(\Gamma)}$ achieves almost the minimum Bayes
risk if $d_1^* / \{\tr(D)-\sum_{j=1}^p d_j^*\} \approx0$.
In terms of Bayes risk reduction, the bound (\ref{bayes-close}) shows that
\begin{eqnarray*}
\tr(D) - R\{\delta_{A^\dag(\Gamma)}, \pi_\Gamma\} \ge\biggl( 1-
\frac{d_1^*+d_2^*+d_3^*+d_4^*}{ \sum_{j=1}^p d_j^*} \biggr) \bigl\{
\tr(D) - R\bigl(\delta^{\mathrm{Bayes}}_\Gamma,
\pi_\Gamma\bigr) \bigr\}.
\end{eqnarray*}
Therefore, $\delta_{A^\dag(\Gamma)}$ achieves Bayes risk reduction
within a negligible factor of that achieved by the Bayes rule if $d_1^*
/ \sum_{j=1}^p d_j^* \approx0$.

In the homoscedastic case where both $D=\sigma^2 I$ and $\Gamma
=\gamma I$, $\delta_{A^\dag}$ reduces to $\delta_{p-2}^{\mathrm{JS}}$,
regardless of $\gamma\ge0$ (Section~\ref{sec3.3}). Then the bounds (\ref
{tight-bound1}) and (\ref{tight-bound2}) become exact and give Efron
and Morris's \cite{EfrMor73} result that $R(\delta_{p-2}^{\mathrm{JS}} , \pi_{\gamma I})
= \tr(D) - (p-2) \{\sigma^4/(\sigma^2+\gamma)\}$ or equivalently $\tr(D) -
R(\delta_{p-2}^{\mathrm{JS}}, \pi_{\gamma I}) = (1-2/p) \{\tr(D) - R(\delta
^{\mathrm{Bayes}}_{\gamma I}, \pi_{\gamma I})\}$.

It is interesting to compare the Bayes risk bound of $\delta_{A^\dag
(\Gamma)}$ with that of the following simpler version of Berger's
\cite{Ber82} estimator $\delta^{\mathrm{MB}}$:
\begin{eqnarray*}
\delta^{\mathrm{MB}2}_j  =
X_j - \Biggl\{ \frac{1}{d_j^*} \sum_{k=j}^p
\bigl(d_k^*-d_{k+1}^*\bigr) \frac{(k-2)_+}{\sum_{\ell=1}^k X_\ell
^2/(d_\ell
+\gamma_\ell)} \Biggr\}
\frac{d_j}{d_j+\gamma_j} X_j.
\end{eqnarray*}
%
By Berger \cite{Ber82}, $\delta^{\mathrm{MB}2}$ is minimax and
%
\begin{eqnarray}
\label{MB-bayes1}R\bigl(\delta^{\mathrm{MB}2}, \pi_\Gamma\bigr) &
=& \tr(D) - \sum_{j=3}^p d_j^*
- 2 \sum_{j=3}^p \frac{d_j^*}{j}
\Biggl( 1 - \frac{d_j^*}{j-1} \sum_{k=1}^{j-1}
\frac{1}{d_k^*} \Biggr)
\\
\label{MB-bayes2}& \le&\tr(D) - \sum_{j=3}^p
d_j^*.
\end{eqnarray}
There seems to be no definite comparison between the bounds (\ref
{tight-bound1}) and (\ref{tight-bound2}) on $R\{\delta_{A^\dag
(\Gamma)}, \pi_\Gamma\}$ and the exact expression (\ref{MB-bayes1})
for $R(\delta^{\mathrm{MB2}}, \pi_\Gamma)$, although the simple bounds (\ref
{loose-bound1}) and (\ref{loose-bound2}) is slightly higher, by at
most $d_3^*+d_4^*$, than the bound (\ref{MB-bayes2}).
Of course, each risk upper bound gives a conservative estimate of the
actual performance, and comparison of two upper bounds should be
interpreted with caution.
In fact, the positive-part estimator $\delta_{A^\dag}^+$ yields lower
risks than those of the non-simplified estimator $\delta^{\mathrm{MB}}$ in our
simulation study (Section~\ref{sec4}).

The simplicity of $\delta_{A^\dag}$ and $\delta_{A^\dag}^+$ makes
it easy to further study them in other ways than using the Bayes (or
average) risk. No similar result to the following Theorem~\ref{th4} has been
established for $\delta^{\mathrm{MB}}$ or $\delta^{\mathrm{MB2}}$.
Corresponding to the prior $\N(0,\Gamma)$, consider the worst-case
(or maximum) risk
\[
R(\delta, \mathcal H_\Gamma) = \sup_{\theta\in\mathcal H_\Gamma}
R(\delta,
\theta)
\]
over the hyper-rectangle $\mathcal H_\Gamma= \{\theta\dvt  \theta_j^2
\le\gamma_j, j=1,\ldots, p\}$ (e.g., Donoho \textit{et al.} \cite{DonLiuMac90}). Applying
Jensen's inequality to (\ref{point-bound}) shows that if $c^*(D,A)>0$, then
\begin{eqnarray*}
R(\delta_A, \theta)  \le\tr(D) - \frac{{c^*}^2(D,A)}{  \sum
_{j=1}^p (d_j+\theta_j^2) a_j^2 } ,
\end{eqnarray*}
which immediately leads to
%
\begin{eqnarray}\label{minimax-bound}
R(\delta_A, \mathcal H_\Gamma)  \le\tr(D) -
\frac{{c^*}^2(D,A)}{
 \sum_{j=1}^p (d_j+\gamma_j) a_j^2 } .
\end{eqnarray}
By the discussion after (\ref{bayes-bound2}), a direct application of
Jensen's inequality to (\ref{bayes-bound}) shows that the Bayes risk
$R(\delta_A, \pi_\Gamma)$ is also no greater than the right-hand
side of (\ref{minimax-bound}), whereas inequality (\ref
{bayes-bound2}) leads to a strictly tighter bound (\ref
{bayes-bound3}). Nevertheless, the upper bound (\ref{minimax-bound})
on the worst-case risk of $\delta_{A^\dag(\Gamma)}$ gives
\begin{eqnarray*}
R \{ \delta_{A^\dag(\Gamma)}, \mathcal H_\Gamma\} \le\tr(D) -
M_\nu= \tr(D) - \Biggl\{ \frac{(\nu-2)^2}{ \sum_{j=1}^\nu
{d_j^*}^{-1}} + \sum
_{j=\nu+1}^p d^*_j \Biggr\},
\end{eqnarray*}
similarly as how (\ref{bayes-bound3}) leads to (\ref{bayes-bound4})
on the Bayes risk of $\delta_{A^\dag(\Gamma)}$. Therefore, the
following result holds by the same proof of Theorem~3.
%
\begin{Theorem}\label{th4} Suppose that $\mathcal H_\Gamma= \{\theta\dvt
\theta_j^2 \le\gamma_j, j=1,\ldots, p\}$. If $\nu=3$, then
\begin{eqnarray*}
R\{\delta_{A^\dag(\Gamma)}, \mathcal H_\Gamma\}  \le\tr(D) -
\sum
_{j=3}^p d_j^* +
\frac{2}{3} d_3^*.
\end{eqnarray*}
If $\nu\ge4$, then
\begin{eqnarray*}
R\{\delta_{A^\dag(\Gamma)}, \mathcal H_\Gamma\} & \le&\tr(D) -
\sum
_{j=3}^p d_j^* + \biggl(
d_3^* + d_4^* - 4\frac{d_\nu^*}{\nu} \biggr)
\\
& \le&\tr(D) - \sum_{j=3}^p
d_j^* + \bigl( d_3^* + d_4^* \bigr).
\end{eqnarray*}
\end{Theorem}

There are similar implications of Theorem~\ref{th4} to those of Theorem~\ref{th3}.
By Donoho \textit{et al.} \cite{DonLiuMac90}, the minimax linear risk over $\mathcal
H_\Gamma$, $R^L(\mathcal H_\Gamma) = \inf_{\delta\,\mathrm{linear}}
R(\delta, \mathcal H_\Gamma)$, coincides with the minimum Bayes risk
$R(\delta^{\mathrm{Bayes}}_\Gamma, \pi_\Gamma)$, and is no greater than
$1.25$ times the minimax risk over $\mathcal H_\Gamma$, $R^N(\mathcal
H_\Gamma)= \inf_{\delta} R(\delta, \mathcal H_\Gamma)$. These
results are originally obtained in the homoscedastic case ($d_1=\cdots
=d_p$), but they remain valid in the heteroscedastic case by the
independence of the observations $X_j$ and the separate constraints on
$\theta_j$. Therefore, a similar result to (\ref{bayes-close}) holds:
\begin{eqnarray*}
R\{ \delta_{A^\dag(\Gamma)}, \mathcal H_\Gamma\} &\le& R^L(
\mathcal H_\Gamma) + \bigl(d_1^*+d_2^*+d_3^*+d_4^*
\bigr)
\\
&\le&1.25 R^N(\mathcal H_\Gamma) + \bigl(d_1^*+d_2^*+d_3^*+d_4^*
\bigr).
\end{eqnarray*}
If $d_1^* / \{\tr(D)-\sum_{j=1}^p d_j^*\} \approx0$, then $\delta
_{A^\dag}$ achieves almost the minimax linear risk (or the minimax
risk up to a factor of $1.25$) over the hyper-rectangle $\mathcal
H_\Gamma$, in addition to being globally minimax with $\theta$ unrestricted.


The foregoing results might be considered non-adaptive in that $\delta
_{A^\dag(\Gamma)}$ is evaluated with respect to the prior $\N
(0,\Gamma)$ or the parameter set $\mathcal H_\Gamma$ with the same
$\Gamma$ used to construct $\delta_{A^\dag(\Gamma)}$. But, by the
invariance of $\delta_A$ under scale transformations of $A$, $\delta
_{A^\dag(\Gamma)}$ is identical to the estimator, $\delta_{A^\dag
(\Gamma_\alpha)}$, that would be obtained if $\Gamma$ is replaced by
$\Gamma_\alpha= \alpha(D+\Gamma)-D$ for any scalar $\alpha$ such
that the diagonal matrix $\Gamma_\alpha$ is nonnegative definite. By
Theorems 3--4, this observation leads directly to the following
adaptive result. In contrast, no adaptive result seems possible for
$\delta^{\mathrm{MB}}$.
%
\begin{Corollary}\label{cor4} 
Let $\Gamma_\alpha= \alpha(D+\Gamma)-D$ and $\alpha_0 = \max
_{j=1,\ldots,p} \{ d_j/(d_j+\gamma_j)\}\ ( \le1)$. Then for each
$\alpha\ge\alpha_0$,
\begin{eqnarray*}
\max\bigl[ R\{ \delta_{A^\dag(\Gamma)}, \pi_{\Gamma_\alpha} \},
R\{
\delta_{A^\dag(\Gamma)}, \mathcal H_{\Gamma_\alpha}\} \bigr] &\le&
R\bigl(
\delta^{\mathrm{Bayes}}_{\Gamma_\alpha}, \pi_{\Gamma_\alpha}\bigr) +
\alpha^{-1} \bigl(d_1^*+d_2^*+d_3^*+d_4^*
\bigr)
\\
&=& R^L(\mathcal H_{\Gamma_\alpha})+ \alpha^{-1} \bigl(
d_1^*+d_2^*+d_3^*+d_4^*\bigr),
\end{eqnarray*}
where $R( \delta^{\mathrm{Bayes}}_{\Gamma_\alpha}, \pi_{\Gamma_\alpha}) =
\tr(D) - \alpha^{-1} \sum_{j=1}^p d_j^*$.
\end{Corollary}

For fixed $\Gamma$, $\delta_{A^\dag(\Gamma)}$ can achieve close to
the minimum Bayes risk or the minimax linear risk with respect to each
prior in the class $\{\N(0,\Gamma_\alpha)\dvt  \alpha\ge\alpha_0\}$
or each parameter set in the class $\{\mathcal H_{\Gamma_\alpha}\dvt
\alpha\ge\alpha_0\}$ under mild conditions. For illustration,
consider the case of a heteroscedastic prior with $\Gamma\propto D$.
Then $\{\Gamma_\alpha\dvt  \alpha\ge\alpha_0\}$ can be reparameterized
as $\{\gamma D\dvt  \gamma\ge0\}$. By Corollary~\ref{cor4}, for each $\gamma\ge0$,
\begin{eqnarray*}
\max\bigl\{R( \delta_{A^\dag_0}, \pi_{\gamma D} ), R(
\delta_{A^\dag
_0}, \mathcal H_{\gamma D})\bigr\} \le R\bigl(
\delta^{\mathrm{Bayes}}_{\gamma D}, \pi_{\gamma D}\bigr) +
\frac{d_1+d_2+d_3+d_4}{1+\gamma} ,
\end{eqnarray*}
where $R( \delta^{\mathrm{Bayes}}_{\gamma D}, \pi_{\gamma D}) = \{\gamma
/(1+\gamma)\} \tr(D)$ and $d_1 \ge d_2 \ge\cdots\ge d_p$.
Therefore, if $d_1/\tr(D) \approx0$, then $\delta_{A^\dag_0}$
achieves the minimum Bayes risk, within a negligible factor, under the
prior $\N(0, \gamma D)$ for each $\gamma>0$. This can be seen as an
extension of the result that in the homoscedastic case, $\delta
_{p-2}^{\mathrm{JS}}$ asymptotically achieves the minimum Bayes risk under the
prior $\N(0, \gamma I)$ for each $\gamma>0$ as $p\to\infty$.

Finally, we compare the estimator $\delta_{A^\dag}$ with a block
shrinkage estimator, suggested by the differentiation in the shrinkage
of low- and high-variance coordinates by $\delta_{A^\dag}$. Consider
the estimator
\begin{eqnarray*}
\delta^{\mathrm{block}} = \lleft\{ %
\begin{array} {c}
\delta_{\tau-2}^{\mathrm{B}}(X_1,\ldots,X_\tau)
\\\noalign{\vspace*{2pt}}
\delta_{p-\tau-2}^{\mathrm{B}} (X_{\tau+1},\ldots,X_p)
\end{array}
\rright\},
\end{eqnarray*}
where $\tau$ is a cutoff index, and $\delta_c^{\mathrm{B}}(Y)=Y$ if $Y$ is of
dimension 1 or 2.
The index $\tau$ can be selected such that the coordinate variances
are relatively homogeneous in each block.
Alternatively, a specific strategy for selecting $\tau$ is to minimize
an upper bound on the Bayes risk of $\delta^{\mathrm{block}}$, similarly as in
the development of $\delta_{A^\dag}$. Applying (\ref{bayes-bound3})
with $A=D^{-1}$ to $\delta_{p-2}^{\mathrm{B}}$ in the two blocks shows that
$R(\delta^{\mathrm{block}}, \pi_\Gamma) \le\tr(D) - L_\tau$, where
\begin{eqnarray*}
L_k = \frac{k-2}{\sklfrac{1}{k}  \sum_{j=1}^k \vfrac{d_j+\gamma
_j}{d_j^2}} + \frac{p-k-2}{({1}/({p-k}))  \sum_{j=k+1}^p \vfrac
{d_j+\gamma_j}{d_j^2}} .
\end{eqnarray*}
The first (or second) term in $L_k$ is set to 0 if $k\le2$ (or $k\ge
p-2$). Then $\tau$ can be defined as the smallest index such that
$L_\tau=\max(L_1,L_2,\ldots,L_p)$. But the upper bound (\ref
{bayes-bound4}) on $R(\delta_{A^\dag}, \pi_\Gamma)$ is likely to be
smaller than the corresponding bound on $R(\delta^{\mathrm{block}}, \pi_\Gamma
)$, because $\{k/(k-2)\} M_k \ge L_k$ for each $k\ge3$ by the
Cauchy--Schwarz inequality $ \{\sum_{j=k+1}^p d_j^2/(d_j+\gamma_j)\} \{
\sum_{j=k+1}^p (d_j+\gamma_j)/d_j^2 \}\ge(p-k)^2$.
Therefore, $\delta_{A^\dag}$ tends to yield greater risk reduction
than $\delta^{\mathrm{block}}$.
This analysis also indicates that $\delta_{A^\dag}$ can be
advantageous over $\delta^{\mathrm{block}}$ extended to multiple blocks.

The rationale of forming blocks in $\delta_{A^\dag}$ and $\delta
^{\mathrm{block}}$ differs from that in existing block shrinkage estimators
(e.g., Brown and Zhao \cite{BroZha09}). As discussed in Cai \cite
{Cai12}, block shrinkage
has been developed mainly in the homoscedastic case as a technique for
pooling information: the coordinate means are likely to be similar to
each other within a block. Nevertheless, it is possible to both deal
with heterogeneity among coordinate variances and exploit homogeneity
among coordinate means within individual blocks in our approach using a
block-homoscedastic prior (i.e., the prior variances are equal within
each block). This topic can be pursued in future work.

\section{Simulation study}\label{sec4}

\subsection{Setup}\label{sec4.1}

We conduct a simulation study to compare the following 8 estimators,
\begin{enumerate}[(ii)]
\item[(i)] Non-minimax estimators: $\delta^{\mathrm{EB}}$ by (\ref{EB}),
$\delta^{\mathrm{XKB}}$ by (\ref{XKB}), $\delta^{\mathrm{RB}}$ by (\ref
{RB}) with $\Gamma=0$;

\item[(ii)] Minimax estimators: $\delta_{p-2}^{\mathrm{B}+}$ by (\ref{B+}),
$\delta^{\mathrm{MB}}$ by (\ref{MB}) with $\Gamma=0$ or $\gamma I$ for some
large $\gamma$, $\delta^+_A$ by (\ref{A+}) with $A=A^\dag_0$ and
$A^\dag_\infty$.
\end{enumerate}
Recall that $A_0^\dag$ corresponds to $\Gamma=0$ or $\Gamma\propto
D$ and $A_\infty^\dag$ corresponds to $\Gamma= \gamma I$ with
$\gamma\to\infty$. In contrast, letting the diagonal elements of
$\Gamma$ tend to $\infty$ in any direction in $\delta^{\mathrm{RB}}$ and
$\delta^{\mathrm{MB}}$ leads to $\delta_0=X$.
Setting $\Gamma$ to 0 or $\infty$ is used here to specify the
relevant estimators, rather than to restrict the prior on $\theta$.

For completeness, we also study the following estimators: $\delta
^{\mathrm{B}+}_{2(p-2)}$ by (\ref{B+}), $\delta^{\mathrm{RB}}$ with $p-2$ replaced by
$2(p-2)$ in (\ref{RB}), $\delta^{\mathrm{MB}}$ with $(k-2)_+$ replaced by
$2(k-2)_+$ in (\ref{MB}), and $\delta^+_A$ with $c^*(D,A)$ replaced
by $2 c^*(D,A)$ in (\ref{A+}), referred to as the alternative versions
of $\delta^{\mathrm{B}+}_{p-2}$, $\delta^{\mathrm{RB}}$, $\delta^{\mathrm{MB}}$, and $\delta
^+_A$ respectively.
The usual choices of the factors, $p-2$, $(k-2)_+$, and $c^*(D,A)$, are
motivated to minimize the risks of the non-positive-part estimators,
but may not be the most desirable for the positive-part estimators. As
seen below, the alternative choices $2(p-2)$, $2(k-2)_+$, and
$2c^*(D,A)$ can lead to risk curves for the positive-part estimators
rather different from those based on the usual choices $(p-2)$,
$(k-2)_+$, and $c^*(D,A)$. Therefore, we compare the estimators $\delta
^{\mathrm{B}+}_{p-2}$, $\delta^{\mathrm{RB}}$, $\delta^{\mathrm{MB}}$, and $\delta^+_A$ and,
separately, their alternative versions.

Each estimator $\delta$ is evaluated by the pointwise risk function
$R(\delta,\theta)$ as $\theta$ moves in a certain direction or the
Bayes risk function $R(\delta,\pi)$ as $\pi$ varies in a set of
priors on $\theta$. Consider the homoscedastic prior $\N(0,\eta^2
I/p)$ or the heteroscedastic prior $\N\{0, \eta^2 D/\tr(D)\}$ for
$\eta\ge0$. As discussed in Section~\ref{sec3.3}, the Bayes risk with the
first or second prior is meant to measure average risk reduction over
the region $\{\theta\dvt  \|\theta\|^2 \le\eta^2\}$ or $\{\theta\dvt
\theta^\T D^{-1} \theta\le p \eta^2/\tr(D)\}$. Corresponding to the
two priors, consider the direction along $(\eta/\sqrt{p},\ldots,\eta
/\sqrt{p})$ or
$(\eta\sqrt{d_1}, \ldots,\eta\sqrt{d_p})/\sqrt{\tr(D)}$, where
$\eta$ gives the Euclidean distance from 0 to the point indexed by
$\eta$.
The two directions are referred to as the homoscedastic and
heteroscedastic directions.
%

We investigate several configurations for $D$, including (\ref
{example}) and
%
\begin{eqnarray}
\label
{example-group3}(d_1,d_2,\ldots,d_{10}) &=&
(40,20,10,5,5,5,1,1,1,1) \quad \mbox{or}
\\
\label{example-group22}& =& (40,20,10,7,6,5,4,3,2,1) \quad \mbox{or}
\\
& =& 5\%, 15\%, \ldots, 95\% \mbox{ quantiles of }
8/\chi^2_3 \mbox{ or } 24/\chi^2_5,
\nonumber
\end{eqnarray}
where $\chi^2_k$ is a chi-squared variable with $k$ degrees of
freedom. In the last case, $(d_1,\ldots,d_{10})$ can be considered a
typical sample from a scaled inverse chi-squared distribution, which is
the conjugate distribution for normal variances. In the case (\ref
{example-group3}), the coordinates may be segmented intuitively into
three groups with relatively homogeneous variances. In the case (\ref
{example-group22}), there is no clear intuition about how the
coordinates should be segmented into groups.

For fixed $D$, the pointwise risk $R(\delta,\theta)$ is computed by
repeatedly drawing $X\sim\N(\theta,D)$ and then taking the average
of $\|\delta-\theta\|^2$. The Bayes risk is computed by repeatedly
drawing $\theta\sim\N(0,\Gamma)$ and $X|\theta\sim\N(\theta,D)$
and then taking the average of $\|\delta-\theta\|^2$. Each Monte
Carlo sample size is set to $10^5$.
%
\begin{figure}

\includegraphics{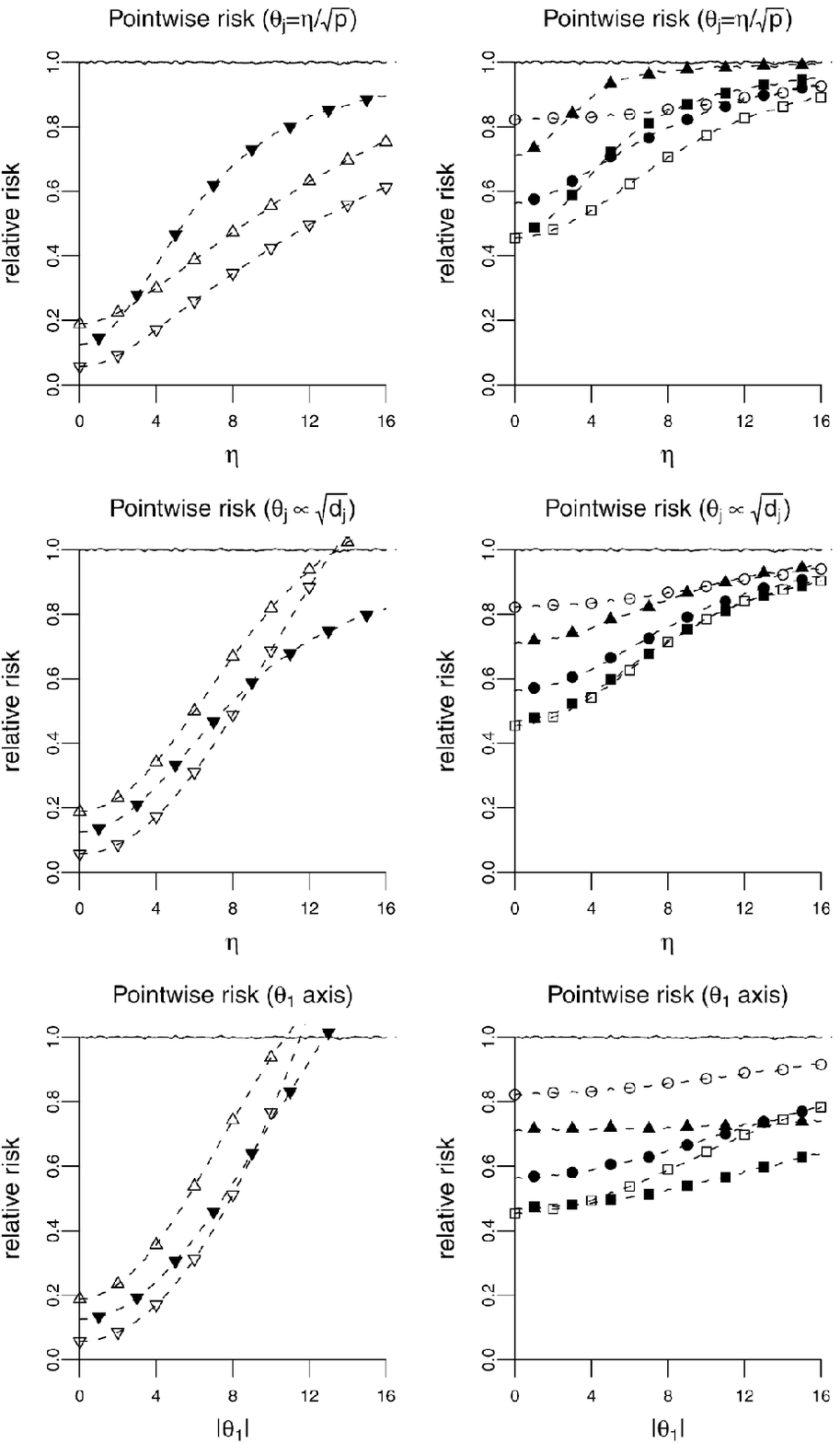}

\caption{Pointwise risks along the homoscedastic (first row) and
heteroscedastic (second row) directions and $\theta_1$ axis (third
row) in the case (\protect\ref{example-group3}). Left: non-minimax estimators
$\delta^{\mathrm{EB}}$ ($\triangledown$), $\delta^{\mathrm{RB}}$ ($\blacktriangledown
$), $\delta^{\mathrm{XKB}}$ ($\vartriangle$). Right: minimax
estimators $\delta_{p-2}^{\mathrm{B}+}$ ($\blacktriangle$), $\delta^{\mathrm{MB}}$
with $\Gamma=0$ ($\bullet$) and $\Gamma=(16^2/p)I$ ($\circ$),
$\delta^+_A$ with $A=A^\dag_0$ ($\blacksquare$) and $A^\dag_\infty
$ ($\square$).}\label{fig1}
\end{figure}
%
%
\begin{figure}

\includegraphics{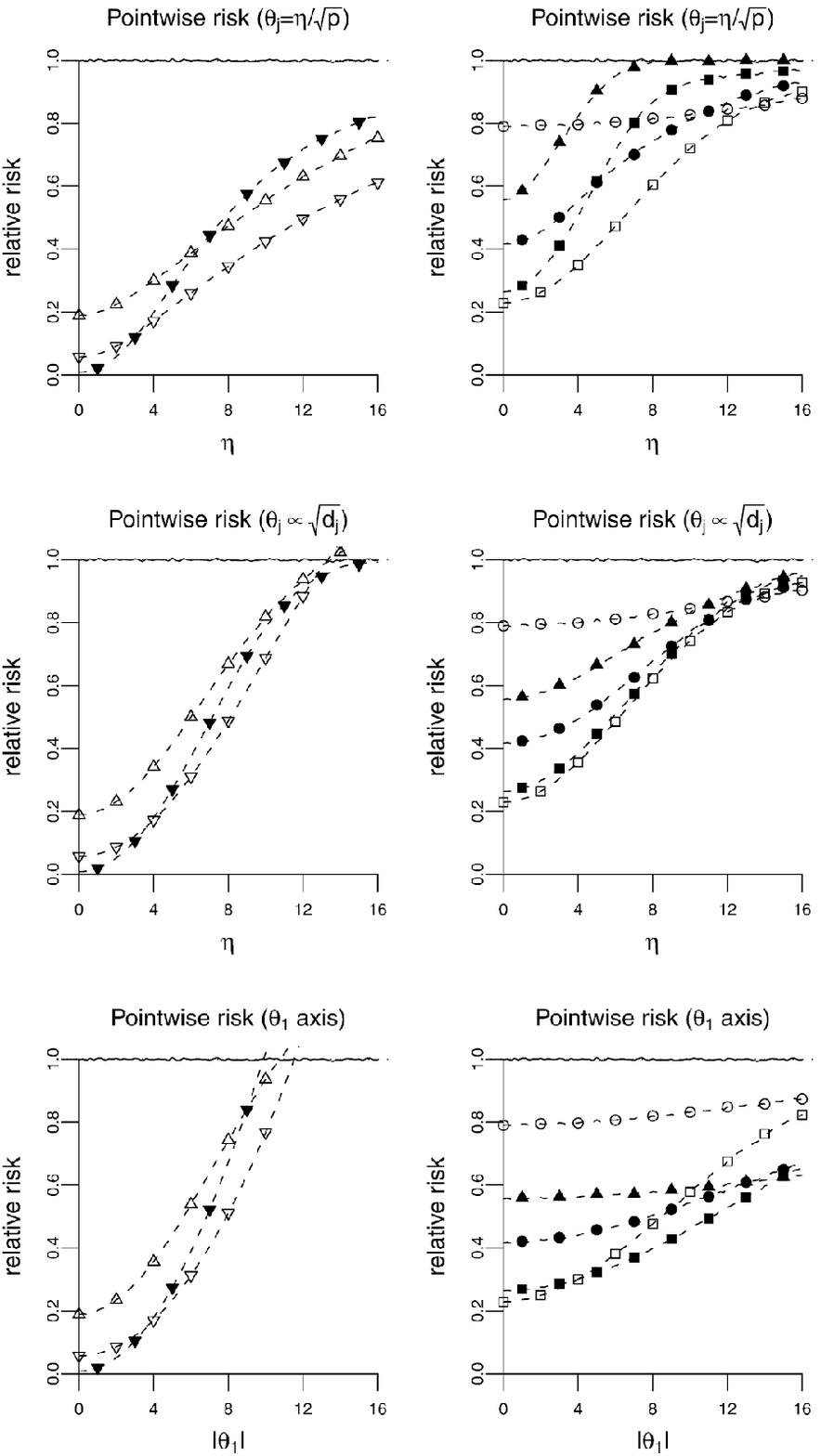}

\caption{Pointwise risks along the homoscedastic (first row) and
heteroscedastic (second row) directions and $\theta_1$ axis (third
row) in the case (\protect\ref{example-group3}), with the same legend
as in
Figure~\protect\ref{fig1}. The alternative versions of $\delta^{\mathrm{B}+}_{p-2}$, $\delta
^{\mathrm{RB}}$, $\delta^{\mathrm{MB}}$, and $\delta^+_A$ are used.}\label{fig2}
\end{figure}

\subsection{Results}

The relative performances of the estimators are found to be consistent
across different configurations of $D$ studied.
Moreover, the Bayes risk curves under the homoscedastic prior are
similar to the pointwise risk curves along the homoscedastic direction.
The Bayes risk curves under the heteroscedastic prior are similar to
the pointwise risk curves along the heteroscedastic direction. Figure~\ref{fig1}
shows the pointwise risks of the estimators with the usual versions of
$\delta^{\mathrm{B}+}_{p-2}$, $\delta^{\mathrm{RB}}$, $\delta^{\mathrm{MB}}$, and $\delta^+_A$
and Figure~\ref{fig2} shows those of the estimators with the alternative
versions of $\delta^{\mathrm{B}+}_{p-2}$, $\delta^{\mathrm{RB}}$, $\delta^{\mathrm{MB}}$, and
$\delta^+_A$ for the case (\ref{example-group3}), with roughly three
groups of coordinate variances, which might be considered unfavorable
to our approach. For both $A^\dag_0$ and $A^\dag_\infty$, the cutoff
index $\nu$ is found to be 3. See the supplementary material (Tan
\cite{Tan}) for the Bayes risk curves of all these estimators for the case
(\ref{example-group3}) and the results for other configurations of $D$.

A number of observations can be drawn from Figures~\ref{fig1}--\ref{fig2}. First, $\delta
^{\mathrm{EB}}$, $\delta^{\mathrm{XKB}}$, and $\delta^{\mathrm{RB}}$ have among the
lowest risk curves along the homoscedastic direction. But along the
heteroscedastic direction, the risk curves of $\delta^{\mathrm{EB}}$ and
$\delta^{\mathrm{XKB}}$ rise quickly above the constant risk of $X$
as $\eta$ increases. Moreover, all the risk curves of $\delta^{\mathrm{EB}}$,
$\delta^{\mathrm{XKB}}$, and $\delta^{\mathrm{RB}}$ along the $\theta_1$
axis exceed the constant risk of $X$ as $|\theta_1|$ increases.
Therefore, $\delta^{\mathrm{EB}}$, $\delta^{\mathrm{XKB}}$, and $\delta
^{\mathrm{RB}}$ fail to be minimax, as mentioned in Section~\ref{sec2}.

Second, $\delta_{p-2}^{\mathrm{B}+}$ or $\delta_{2(p-2)}^{\mathrm{B}+}$ has among the
highest risk curve, except where the risk curves of $\delta^{\mathrm{EB}}$ and
$\delta^{\mathrm{XKB}}$ exceed the constant risk of $X$ along the
heteroscedastic direction. The poor performance is expected for $\delta
_{p-2}^{\mathrm{B}+}$ or $\delta_{2(p-2)}^{\mathrm{B}+}$, because there are considerable
differences between the coordinate variances in (\ref{example-group3}).

Third, among the minimax estimators, $\delta^+_A$ with $A=A^\dag_0$
or $A^\dag_\infty$ has the lowest risk curve along various
directions, whether the usual versions of $\delta^{\mathrm{B}+}_{p-2}$, $\delta
^{\mathrm{MB}}$, and $\delta^+_A$ are compared (Figure~\ref{fig1}) or the alternative
versions are compared (Figure~\ref{fig2}).

Fourth, the risk curve of $\delta^+_A$ with $A=A^\dag_0$ is similar
to that of $\delta^+_A$ with $A=A^\dag_\infty$ along the
heteroscedastic direction. But the former is noticeably higher than the
latter along the homoscedastic direction as $\eta$ increases, whereas
is noticeably lower than the latter along the $\theta_1$ axis as
$|\theta_1|$ increases.
These results agree with the construction of $A^\dag_0$ using a
heteroscedastic prior and $A^\dag_\infty$ using a flat, homoscedastic
prior. Their relative performances depend on the direction in which the
risks are evaluated.

Fifth, $\delta^{\mathrm{MB}}$ with $\Gamma=0$ has risk curves below that of
$\delta_{p-2}^{\mathrm{B}+}$ or $\delta_{2(p-2)}^{\mathrm{B}+}$, but either above or
crossing those of $\delta^+_A$ with $A=A^\dag_0$ and $A^\dag_\infty
$. Moreover, $\delta^{\mathrm{MB}}$ with $\Gamma=(16^2/p)I$ has elevated,
almost flat risk curves for $\eta$ from 0 to 16. This seems to
indicate an undesirable consequence of using a non-degenerate prior for
$\delta^{\mathrm{MB}}$ in that the risk tends to increase for $\theta$ near 0,
and remains high for $\theta$ far away from 0.

The foregoing discussion involves the comparison of the risk curves as
$\theta$ moves away from 0 between $\delta^{\mathrm{MB}}$ and $\delta_{A^\dag
}^+$ specified with fixed priors. Alternatively, we compare the
pointwise risks at $\theta=(\eta/\sqrt{p},\ldots,\eta/\sqrt{p})$
or $(\eta\sqrt{d_1}, \ldots,\eta\sqrt{d_p})/\sqrt{\tr(D)}$
and the Bayes risks under the prior $\N(0,\eta^2 I/p)$ or $\N\{0,
\eta^2 D/\tr(D)\}$ between $\delta^{\mathrm{MB}}$ and $\delta_{A^\dag}^+$
specified with the prior $\N(0,\eta^2 I/p)$ for a range of $\eta$.
The homoscedastic prior used in the specification of $\delta^{\mathrm{MB}}$ and
$\delta_{A^\dag}^+$ can be considered correctly specified or
misspecified, when the Bayes risks are evaluated under, respectively,
the homoscedastic or heteroscedastic prior or when the pointwise risks
are evaluated along the homoscedastic or heteroscedastic direction. For
each situation, $\delta_{A^\dag}^+$ has lower pointwise or Bayes
risks than $\delta^{\mathrm{MB}}$. See Figure A2 in the supplementary material
(Tan \cite{Tan}).

\section{Conclusion}\label{sec5}

The estimator $\delta_{A^\dag}$ and its positive-part version $\delta
^+_{A^\dag}$ are not only minimax and but also have desirable
properties including simplicity, interpretability, and effectiveness in
risk reduction. In fact, $\delta_{A^\dag}$ is defined by taking
$A=A^\dag$ in a class of minimax estimators $\delta_A$. The
simplicity of $\delta_{A^\dag}$ holds because $\delta_A$ is of the
linear form $(I-\lambda A)X$, with $A$ and $\lambda$ indicating the
direction and magnitude of shrinkage. The interpretability of $\delta
_{A^\dag}$ holds because the form of $A^\dag$ indicates that one
group of coordinates are shrunk in the direction of Berger's \cite{Ber76}
minimax estimator whereas the remaining coordinates are shrunk in the
direction of the Bayes rule. The effectiveness of $\delta_{A^\dag}$
in risk reduction is supported, in theory, by showing that $\delta
_{A^\dag}$ can achieve close to the minimum Bayes risk simultaneously
over a scale class of normal priors (Corollary~\ref{cor4}).
For various scenarios in our numerical study, the estimators $\delta
_{A^\dag}^+$ with extreme priors yield more substantial risk reduction
than existing minimax estimators.

It is interesting to discuss a special feature of $\delta_{A,r}$ and
hence of $\delta_{A,c}$ and $\delta_A$ among linear, shrinkage
estimators of the form
%
\begin{eqnarray}\label{general}
\delta= X - h\bigl(X^\T B X\bigr) A X,
\end{eqnarray}
where $A$ and $B$ are nonnegative definite matrices and $h(\cdot)$ is
a scalar function.
The estimator $\delta_{A,r}$ corresponds to the choice $B \propto A^\T
Q A$, which is motivated by the form of the optimal $\lambda$ in
minimizing the risk of $(I-\lambda A) X$ for fixed $A$. On the other
hand, Berger and Srinivasan \cite{BerSri78} showed that under certain regularity
conditions on $h(\cdot)$, an estimator (\ref{general}) can be
generalized Bayes or admissible only if $B \propto\Sigma^{-1} A$.
This condition is incompatible with $B \propto A^\T Q A$, unless $A
\propto Q^{-1} \Sigma^{-1}$ as in Berger's \cite{Ber76} estimator.
Therefore, $\delta_A$ including $\delta_{A^\dag}$ is, in general,
not generalized Bayes or admissible. This conclusion, however, does not
apply directly to the positive-part estimator $\delta_A^+$, which is
no longer of the linear form $(I-\lambda A)X$.

There are various topics that can be further studied.
First, the prior on $\theta$ is fixed, independently of data in the
current paper. A useful extension is to allow the prior to be
estimated within a certain class, for example, homoscedastic priors
$N(0,\gamma I)$, from the data, in the spirit of empirical Bayes
estimation (e.g., Efron and Morris \cite{EfrMor73}).
Second, the Bayes risk with a normal prior is used to measure average
risk reduction in an elliptical region (Section~\ref{sec3.3}). It is interesting
to study how our approach can be extended when using a non-normal prior
on $\theta$, corresponding to a non-elliptical region in which risk
reduction is desired.
\begin{appendix}
\section*{Appendix}\label{app}

\begin{Preparation*} The following extends Stein's \cite{Ste81} lemma for
computing the expectation of the inner product of $X-\theta$ and a
vector of functions of $X$.
\end{Preparation*}
%
\begin{Lemma}\label{lem1} Let $X=(X_1,\ldots,X_p)^\mathrm{T}$ be multivariate normal
with mean $\theta$ and variance matrix $\Sigma$.
Assume that $g=(g_1,\ldots,g_p)^\T\dvtx  \mathcal R^p \to\mathcal R^p$ is
almost differentiable Stein \cite{Ste81} with $E_\theta
\{ |\nabla_j g_i(X)|
\}<\infty$ for $i,j=1,\ldots,p$, where $\nabla_j=\partial/\partial
x_j$. Then
\[
E_\theta\bigl\{ (X-\theta)^\T g(X) \bigr\} = \tr\bigl[
\Sigma E_\theta\bigl\{\nabla g(X)\bigr\} \bigr],
\]
where $\nabla g(x)$ is the matrix with $(i,j)$th element $\nabla_j g_i(x)$.
\end{Lemma}
\begin{pf} A direct generalization of Lemma~2 in
Stein \cite{Ste81} to a normal random vector with
non-identity variance matrix gives
\[
E_\theta\bigl\{ (X-\theta) g_i(X) \bigr\} = \Sigma
E_\theta^\T\bigl\{ \nabla g_i(X) \bigr\},
\]
where $\nabla g_i(x)$ is the row vector with $j$th element $\nabla_j
g_i(x)$. Taking the $i$th element of both sides of the equation gives
\[
E_\theta\bigl\{ (X_i-\theta_i)
g_i(X) \bigr\} = \sum_{j=1}^p
\sigma_{ij} E_\theta\bigl\{\nabla_j
g_i(X)\bigr\},
\]
where $\sigma_{ij}$ is the $(i,j)$th element of $\Sigma$. Summing
both sides of the preceding equation over $i$ gives the desired result.
\end{pf}
\begin{pf*}{Proof of Theorem~\ref{th1}} By direct calculation, the risk of
$\delta_{A,r}$ is
\begin{eqnarray*}
R(\delta_{A,r},\theta) = \tr(\Sigma Q) + E_\theta\biggl(
\frac
{r^2}{X^\T A^\T Q A X} \biggr) - 2 E_\theta\biggl\{ (X-\theta)^\T
\frac{r Q A X }{X^\T A^\T Q A X} \biggr\}.
\end{eqnarray*}
By Lemma~\ref{lem1} and the fact that $\tr(\Sigma QAX X^\T A^\T QA)= X^\T A^\T
QA\Sigma QAX $, the third term after the minus sign in $R(\delta
_{A,r},\theta)$ is
\begin{eqnarray*}
2E_\theta\biggl\{ r \frac{\tr(\Sigma Q A)}{X^\T A^\T Q A X} \biggr
\} - 4 E_\theta
\biggl\{ r\frac{ X^\T A^\T QA \Sigma QAX }{(X^\T A^\T
Q A X)^2} \biggr\} + 4 E_\theta\biggl(r'
\frac{ X^\T A^\T QA \Sigma
QAX }{X^\T A^\T Q A X} \biggr).
\end{eqnarray*}
By condition (\ref{A-cond}), $A^\T QA \Sigma QA $ is nonnegative
definite. By Section~21.14 and Exercise 21.32 in Harville
\cite{Har08},
$(x^\T A^\T QA \Sigma QA x)/(x^\T A^\T Q A x) \le\lambda_{\max
}(A\Sigma Q + \Sigma A^\T Q)/2$ for $x \neq0$. Then the preceding
expression is bounded from below by
\begin{eqnarray*}
2E_\theta\biggl\{ r \frac{\tr(\Sigma Q A)-\lambda_{\max}(A\Sigma Q
+ \Sigma A^\T Q)}{X^\T A^\T Q A X} \biggr\},
\end{eqnarray*}
which leads immediately to the upper bound on $R(\delta_{A,r},\theta
)$. \end{pf*}
\begin{pf*}{Proof for condition (\ref{A-cond2})} We show that if
condition (\ref{A-cond2}) holds, then there exists a nonsingular
matrix $B$ with the claimed properties. The converse is trivially true.
Let $R$ be the unique symmetric, positive definite matrix such that $R^2=Q$.
Then $R A R^{-1}$ is symmetric, that is, $R A R^{-1} = R^{-1} A^\T R$,
because $QA=A^\T Q$. Moreover, $R \Sigma R$ and $R A R^{-1}$ commute,
that is, $RAR^{-1} (R \Sigma R) =R \Sigma R (RAR^{-1})^\T= R \Sigma R
(RAR^{-1})$, because $A \Sigma=\Sigma A^\T$ and $R A R^{-1}$ is
symmetric. Therefore, $R\Sigma R$ and $R A R^{-1}$ are simultaneously
diagonalizable (Harville \cite{Har08}, Section~21.13).
There exists an
orthogonal matrix $O$ such that $O(R\Sigma R)O^\T=D$ and $O(R A
R^{-1})O^\T= A^*$ for some diagonal matrices $D$ and $A^*$. Then
$B=OR$ satisfies the claimed properties.
\end{pf*}
\begin{pf*}{Proof of inequality (\ref{bayes-bound2})} We show that if
$(Z_1,\ldots,Z_p)$ are independent standard normal variables, then
$E\{ (\sum_{j=1}^p a_j^2 Z_j^2)^{-1}\} \ge\{p/(p-2)\} (\sum_{j=1}^p
a_j^2)^{-1}$. Let $S= \sum_{j=1}^p Z_j^2$. Then $S$ and $(Z_1^2/S,
\ldots, Z_p^2/S)$ are independent, $S\sim\chi^2_p$, and $(Z_1^2/S,
\ldots, Z_p^2/S) \sim\operatorname{Dirichlet}(1/p,\ldots, 1/p)$.
The claimed
inequality follows because $E\{ (\sum_{j=1}^p a_j^2 Z_j^2)^{-1}\} =
E\{ (\sum_{j=1}^p a_j^2 Z_j^2/S)^{-1}\} E(S^{-1})$,
$E(S^{-1})=1/(p-2)$, and
$E\{ (\sum_{j=1}^p$ $ a_j^2 Z_j^2/S)^{-1}\} \ge\break (\sum_{j=1}^p
a_j^2/p)^{-1}$ by Jensen's inequality.
\end{pf*}
\begin{pf*}{Proofs of Theorem~\ref{th2} and Corollary~\ref{cor2}} Consider the
transformation $\delta_j = d_j^2 /(d_j+\gamma_j)$ and $\alpha_j = \{
(d_j+\gamma_j)/d_j\} a_j$, so that $\delta_j \alpha_j = d_j a_j$ and
$\delta_j \alpha_j^2 = (d_j+\gamma_j) a_j^2 $. Problem (\ref{opt})
is then transformed to $
\max_{\alpha_1,\ldots,\alpha_p} \{ \sum_{j=1}^p \delta_j \alpha
_j - 2 \max(\delta_1 \alpha_1, \ldots, \delta_p \alpha_p)\}$,
subject to $\alpha_j\ge0$ ($j=1,\ldots, p$) and $\sum_{j=1}^p
\delta_j \alpha_j^2 = \sum_{j=1}^p \delta_j$, which is of the form
of the special case of (\ref{opt}) with $\gamma_j=0$ ($j=1,\ldots,p$).
But it is easy to verify that if the claimed results hold for the
transformed problem, then the results hold for original problem (\ref
{opt}). Therefore, assume in the rest of proof that $\gamma_j=0$
($j=1,\ldots,p$).

There exists at least a solution, $A^\dag$, to problem (\ref{opt}) by
boundedness of the constraint set. Let $\mathcal K=\{k\dvt  d_k a^\dag_k =
d_\nu a^\dag_\nu, k=1, \ldots, p\}$ and $\mathcal K^c=\{j\dvt  d_j
a^\dag_j < d_\nu a^\dag_\nu, j=1, \ldots, p\}$. A key of the proof
is to exploit the fact that, by the setup of problem (\ref{opt}),
$(a_1^\dag, \ldots, a_p^\dag)$ is automatically a solution to the problem
\setcounter{equation}{0}
\begin{eqnarray}\label{a1}
&&\max_{a_1,\ldots,a_p} \quad  \sum_{j=1}^p
d_j a_j - 2 d_\nu a_\nu,\nonumber
\\
&&\quad \mbox{subject to} \quad  a_j\ge0,\qquad  d_ja_j \le
d_\nu a_\nu\qquad (j=1,\ldots, p), \quad \mbox{and}\\
&&\hphantom{\quad \mbox{subject to} \quad}  \sum
_{j=1}^p d_j a_j^2
= \sum_{j=1}^p d_j.
\nonumber
\end{eqnarray}
The Karush--Kuhn--Tucker condition for this problem gives
%
\begin{eqnarray}
\label{a2}-1 + 2\lambda a_j^\dag-d_j^{-1}
\rho_j &=& 0 \qquad \mbox{for} j \in\mathcal K^c,
\\
\label{a3}-1 + 2\lambda a_k^\dag+ \mu_k-d_k^{-1}
\rho_k &=& 0\qquad  \mbox{for } k\ (\neq\nu) \in\mathcal K,
\\
\label{a4}-1 + 2\lambda a_\nu^\dag+ \biggl(2-\sum
_{k\in\mathcal K \setminus\{
\nu\}} \mu_k \biggr)-d_\nu^{-1}
\rho_\nu& =& 0,
\end{eqnarray}
where $\lambda$, $\mu_k \ge0$ ($k\in\mathcal K\setminus\{\nu\}$),
and $\rho_j\ge0$ satisfying $\rho_j a_j^\dag=0$ ($j=1,\ldots,p$)
are Lagrange multipliers.

First, we show that $a_j^\dag>0$ and hence $\rho_j=0$ for $j=1,\ldots
,p$. If $\mathcal K^c = \emptyset$, then either $a_j^\dag>0$ for
$j=1,\ldots,p$, or $a_1^\dag= \cdots= a_p^\dag=0$. The latter case
is infeasible by the constraint $\sum_{j=1}^p d_j a_j^2 = \sum
_{j=1}^p d_j$. Suppose $\mathcal K^c \neq\emptyset$. By (\ref{a2}),
$a_j^\dag>0$ for each $j\in\mathcal K^c$. Then $a_k^\dag> 0$ for
each $k \in\mathcal K$ because $d_k a_k^\dag> d_j a_j^\dag$.

Second, we show that $\nu\ge3$. If $\mathcal K^c = \emptyset$, then
$\nu=p\ge3$. Suppose $\mathcal K^c \neq\emptyset$. Then $\lambda
>0$ by (\ref{a2}).
Summing (\ref{a3}) over $k\ (\neq\nu) \in\mathcal K $ and (\ref{a4}) shows that
$-\nu+ 2\lambda\sum_{k=1}^\nu a_k^\dag+ 2=0$. Therefore, $\nu>2$
or equivalently $\nu\ge3$.

Third, we show that $\mathcal K=\{1,2,\ldots,\nu\}$ and $\mathcal
K^c=\{\nu+1,\ldots, p\}$. For each $k\ (\neq\nu)\in\mathcal K$ and
$j \in\mathcal K^c$, $a_k^\dag\le a_j^\dag$ by (\ref{a2})--(\ref{a3}) and then
$d_k > d_j$ because $d_k a_k^\dag> d_j a_j^\dag$. The inequalities
also hold for $k=\nu$, by application of the argument to problem (\ref{a1})
with $\nu$ replaced by some $k\ (\neq\nu)\in\mathcal K$. Then
$\mathcal K^c=\{\nu+1,\ldots, p\}$ because $d_\nu>d_j$ for each
$j\in\mathcal K^c$, $d_1\ge d_2 \ge\cdots\ge d_p$, and $\nu$ is the
largest element in $\mathcal K$.

Fourth, we show the expressions for $(a_1^\dag,\ldots,a_p^\dag)$ and
the achieved maximum value.
By the definition of $\mathcal K$, $a_k^\dag\propto d_k^{-1}$ for
$k=1,\ldots,\nu$.
By (\ref{a2}), $a_j^\dag\propto1$ for $j=\nu+1,\ldots,p$.
Let $y^\dag=d_\nu a^\dag_\nu$ and $z^\dag= a^\dag_{\nu+1}$. Then
$(y^\dag,z^\dag)$ is a solution to the problem
\begin{eqnarray*}
&&\max_{y,z} \quad  (\nu-2) y + \Biggl( \sum
_{j=\nu+1}^p d_j \Biggr) z ,
\\
&&\quad \mbox{subject to}\quad   y\ge0,\qquad  z\ge0,\qquad  y \ge d_{\nu+1}z, \quad \mbox{and}\\
&&\hphantom{\quad \mbox{subject to}\quad}\Biggl(
\sum_{k=1}^\nu d_k^{-1}
\Biggr) y^2 + \Biggl( \sum_{j=\nu
+1}^p
d_j \Biggr) z^2 = \sum_{j=1}^p
d_j.
\end{eqnarray*}
By the definition of $\mathcal K$, $y^\dag>d_{\nu+1}z^\dag$ and
hence $(y^\dag, z^\dag)$ lies off the boundary in the constraint set.
Then $(y^\dag,z^\dag)$ is a solution to the foregoing problem with
the constraint $y \ge d_{\nu+1}z$ removed. The problem is of the form
of maximizing a linear function of $(y,z)$ subject to an elliptical constraint.
Straightforward calculation shows that
\begin{eqnarray*}
y^\dag = \biggl( \frac{ \sum_{j=1}^p d_j}{ M_\nu} \biggr)^{1/2}
\frac{\nu-2}{ \sum_{j=1}^\nu d_j^{-1}}, \qquad z^\dag= \biggl( \frac{ \sum
_{j=1}^p d_j}{ M_\nu}
\biggr)^{1/2},
\end{eqnarray*}
and the achieved maximum value is $(\sum_{j=1}^p d_j)^{1/2} M_\nu
^{1/2}$, where
$M_\nu= (\nu-2)^2 /(\sum_{j=1}^\nu d_j^{-1}) $ $+ \sum_{j=\nu+1}^p
d_j $ .

Finally, we show that the sequence $(M_3, M_4, \ldots M_p)$ is
nonincreasing: $M_k \ge M_{k+1}$, where the equality holds if and only
if $k-2 = \sum_{j=1}^k d_{k+1}/d_j$. Because $y^\dag>d_{\nu+1}z^\dag
$ or $\nu-2 > \sum_{j=1}^\nu d_{\nu+1}/d_j$, this result implies
that $M_\nu> M_{\nu+1}$ and hence $A^\dag$ is a unique solution to
(\ref{opt}). Let
$L_k = \{(\sum_{j=1}^k d_j)( \sum_{j=1}^k d_j^{-1}) -(k-2)^2 \}/ \sum
_{j=1}^k d_j^{-1}$ so that $M_k = \sum_{j=1}^p d_j - L_k$.
By the identity $(b+\beta)/(a+\alpha) - b/a = (\beta/ \alpha-
b/a)\{ \alpha/(a+\alpha)\}$ and simple calculation,
%
\begin{eqnarray}\label{a5}
L_{k+1} - L_k & =& \Biggl[ \frac{  \sum_{j=1}^k (\sfrac{d_j}{d_{k+1}} +
\sfrac{d_{k+1}}{d_j}) - 2k +4}{ d_{k+1}^{-1}} - \Biggl\{
\sum_{j=1}^k d_j -
\frac{(k-2)^2 }{  \sum_{j=1}^k d_j^{-1}} \Biggr\} \Biggr] \frac
{d_{k+1}^{-1}}{\sum_{j=1}^{k+1} d_j^{-1}}
\nonumber
\\[-6pt]
\\[-10pt]
& =& d_{k+1} \frac{ \{ r_k - (k-2) \}^2}{r_k (r_k+1)},\nonumber
\end{eqnarray}
where $r_k = \sum_{j=1}^k d_{k+1}/d_j$. Therefore, $L_k \le L_{k+1}$.
Moreover, $L_k=L_{k+1}$ if and only if $r_k=k-2$, that is, $\sum
_{j=1}^k d_{k+1}/d_j=k-2$.
\end{pf*}
\begin{pf*}{Proof of Corollary~\ref{cor3}} It suffices to show (\ref{sol-ineq}).
By Corollary~\ref{cor2}, $\sum_{k=1}^{\nu-1} d_\nu^*/ d_k^* \ge\nu-3$ and hence
$\sum_{k=1}^\nu d_\nu^*/ d_k^* \ge\nu-2$. Then for $j=1,\ldots,\nu$.
\begin{eqnarray*}
a_j^\dag= \frac{(\nu-2) {d_j^*}^{-1} }{  \sum_{k=1}^\nu{d_k^*}^{-1}
} \frac{d_j}{d_j+\gamma_j} \le
\frac{(\nu-2) {d_\nu^*}^{-1} }{
 \sum_{k=1}^\nu{d_k^*}^{-1} } \frac{d_j}{d_j+\gamma_j} \le\frac
{d_j}{d_j+\gamma_j},
\end{eqnarray*}
because $d_j^* \ge d_\nu^*$ for $j \le\nu$.
\end{pf*}
\begin{pf*}{Proof of Theorem~\ref{th3}} Let
$L_k = \sum_{j=1}^k d_j^* -(k-2)^2/ \sum_{j=1}^k {d_j^*}^{-1}$ so
that $M_k = \sum_{j=1}^p d_j^* - L_k$, similarly as in the proof of
Theorem~\ref{th2}. By equation (\ref{a5}) with $r_k = \sum_{j=1}^k d_{k+1}^*/d_j^*$
and $d_{k+1}$ replaced by $d_{k+1}^*$,
\begin{eqnarray*}
L_\nu= L_3 + \sum_{k=3}^{\nu-1}
(L_{k+1} - L_k)= L_3 + \sum
_{k=3}^{\nu-1} d_{k+1}^* \frac{ \{ r_k - (k-2) \}^2}{r_k (r_k+1)} .
\end{eqnarray*}
By the relationship $r_k = (d_{k+1}^* / d_k^*) (1+r_{k-1})$ and simple
calculation,
\begin{eqnarray*}
L_3 & =& d_1^* + d_2^* + d_3^*
- \frac{1}{{d_1^*}^{-1} + {d_2^*}^{-1} +
{d_3^*}^{-1}}
\\
& =& d_1^* + d_2^* + d_3^* - \sum
_{k=3}^{\nu-1} d_{k+1}^* \biggl(
\frac{1}{r_k} - \frac{1}{r_k+1} \biggr) - \frac{d_\nu^*}{r_{\nu-1}+1}.
\end{eqnarray*}
If $\nu\ge4$, combining the two preceding equation gives
\begin{eqnarray*}
L_\nu&=& d_1^* + d_2^* + d_3^* +
\sum_{k=3}^{\nu-1} d_{k+1}^*
\frac{
\{ r_k - (k-2) \}^2-1}{r_k (r_k+1)} - \frac{d_\nu^*}{r_{\nu-1}+1}
\\
& \le &d_1^* + d_2^* + d_3^* + \sum
_{k=3}^{\nu-1} d_{k+1}^*
\frac
{3}{k(k+1)} - \frac{d_\nu^*}{\nu}
\\
& =& d_1^* + d_2^* + d_3^* +
d_4^* - 3 \sum_{k=3}^{\nu-2}
\frac
{d_{k+1}^* - d_{k+2}^*}{k+1} - 4 \frac{d_\nu^*}{\nu}
\\
& \le& d_1^* + d_2^* + d_3^* +
d_4^* -4 \frac{d_\nu^*}{\nu}.
\end{eqnarray*}
The first inequality follows because $k-2 \le r_k \le k$ for
$k=3,\ldots,\nu-1$ and $\{t - (k-2)\}^2 /\{t(t+1)\}$ is increasing
for $k-2 \le t \le k$ with a maximum at $t=k$. The second inequality
follows because $d_1^* \ge d_2^* \ge\cdots\ge d_p^*$. Therefore, if
$\nu\ge4$ then
\begin{eqnarray*}
\frac{p}{p-2} M_\nu& \ge&\frac{p}{p-2} \Biggl\{ \sum
_{j=1}^p d_j^* - \biggl(
d_1^* + d_2^* + d_3^* + d_4^*
-4 \frac{d_\nu^*}{\nu} \biggr) \Biggr\}
\\
& =& \sum_{j=3}^p
d_j^* - \Biggl( d_3^* + d_4^* -
\frac{2}{p-2} \sum_{j=5}^p
d_j^* -\frac{ 4 p}{p-2} \frac{d_\nu^*}{\nu} \Biggr) .
\end{eqnarray*}
If $\nu=3$, then $L_\nu\le d_1^* + d_2^* + d_3^* -d_3^*/3$ and hence
\begin{eqnarray*}
\frac{p}{p-2} M_\nu& \ge&\frac{p}{p-2} \Biggl\{ \sum
_{j=1}^p d_j^* -
\bigl(d_1^* + d_2^* + d_3^*
-d_3^*/3 \bigr) \Biggr\}
\\
& =&\sum_{j=3}^p d_j^* -
\Biggl( d_3^* - \frac{2}{p-2} \sum_{j=4}^p
d_j^* -\frac{ p}{p-2} \frac{d_3^*}{3} \Biggr). 
\end{eqnarray*}
This completes the proof. \end{pf*}
\end{appendix}

%

\section*{Acknowledgements}
The author thanks Bill Strawderman and Cunhui Zhang for helpful discussions.

\begin{supplement}
\stitle{Supplementary Material for ``Improved minimax estimation of a multivariate normal mean under heteroscedasticity''}
\slink[doi]{10.3150/13-BEJ580SUPP} 
\sdatatype{.pdf}
\sfilename{BEJ580\_supp.pdf}
\sdescription{We present additional results from the simulation study in Section~\ref{sec4}.}
\end{supplement}

%

\printhistory

\end{document}